\documentclass{article}
\usepackage{amsmath}
\usepackage[latin1]{inputenc}
\usepackage[T1]{fontenc}
\usepackage{lmodern}
\usepackage{fullpage}
\usepackage[english]{babel}
\usepackage{amsfonts}
\usepackage{graphicx}
\usepackage{verbatim}
\usepackage{color}
\definecolor{liens}{rgb}{0.,0.,0.8}
\usepackage[colorlinks,
            linkcolor=liens,
            anchorcolor=liens,
            citecolor=liens,
            filecolor=liens,
            menucolor=liens,
            urlcolor=liens,
            bookmarks,
            bookmarksnumbered,
            pdfstartview=FitH]{hyperref}
\usepackage{multirow}

\newcommand{\RR}{\mathbb R}
\newcommand{\NN}{\mathbb N}
\newcommand{\DD}{\mathbb D}
\newcommand{\CC}{\mathbb C}
\newcommand{\TT}{\mathbb T}

\newcommand{\ZZ}{\mathbb Z}
\newcommand{\sph}{\mathbb S}

\newcommand{\cP}{{\cal P}}
\newcommand{\cR}{{\cal R}}
\newtheorem{theorem}{Theorem}
\newcommand{\boite}{\mbox{} \hfill \mbox{\rule{2mm}{2mm}}}
\newtheorem{prop}{Proposition}
\newtheorem{lem}{Lemma}

\newtheorem{corollary}{Corollary}

\begin{document}

\title{Minimax principle and lower bounds in  $H^2$-rational approximation\footnote{This work has been partly funded by \emph{Macao Government FDCT 098/2012/A3.}}}
\author{Laurent Baratchart\footnotemark[2] ~~~~~~~~~ Sylvain Chevillard\footnotemark[2] ~~~~~~~~~ Tao Qian\footnotemark[3] \\[0.3cm]
\url{Laurent.Baratchart@sophia.inria.fr},~phone: +33 4 92 38 78 74, fax: +33 4 92 38 78 58.\\
\url{sylvain.chevillard@inria.fr},~phone: +33 4 92 38 76 42, fax: +33 4 92 38 78 58.\\
\url{fsttq@umac.mo},~phone:  +853 88 22 85 47, fax: +853 88 22 24 26.
}

\footnotetext[2]{Inria, 2004 route des Lucioles, BP 93, 06\,902 Sophia-Antipolis Cedex, France.}
\footnotetext[3]{Faculty of Science and Technology, University of Macau, E11, Avenida da Universidade, Taipa, Macau, China.}
\date{}

\maketitle

\paragraph{Abstract:}
We derive  lower bounds in rational approximation of given degree 
to functions in the Hardy space $H^2$ of the unit disk. We apply these to 
asymptotic errors rates in rational
approximation to Blaschke products and to Cauchy integrals on geodesic arcs.
We also explain how to compute such bounds, 
either using Adamjan-Arov-Krein theory or linearized errors, 
and we present a couple of numerical  experiments. 
We dwell on 
a maximin principle 
developed in \cite{BS02}.

{\it Classification numbers (AMS):}  31B05, 35J25, 42B35, 46E20, 47B35.\\

\tableofcontents

\newpage
\section{Introduction}
Rational approximation to a given function on a curve in the complex plane
is a classical topic from analysis, and a cornerstone
of modeling and design in several areas of applied sciences and engineering.
Special interest attaches to the case where the approximated function 
extends holomorphically on one side of the curve.
In connection with system identification and control,
such issues typically arise on the line or the circle where they make 
contact  with extremal problems in
Hardy spaces \cite{B_CMFT99,GL84,Peller,Nikolskii2,OSM,IIS,YuWang}. 
Our model curve in this paper will be the circle, though everything 
translates easily to the line. The criterion under examination 
will be the $L^2$-norm.

From the 
approximation-theoretic viewpoint, much attention has been directed towards
error rates, in connection with  smoothness of the approximated function.
Let us mention Peller's converse theorems 
on the speed of rational approximation
\cite{Peller}, Glover's construction of near-best uniform rational
approximants \cite{GL84}, Parfenov's solution of a
conjecture by Gonchar on the degree of rational approximation to
holomorphic functions on compact subsets of the domain of analyticity
\cite{Par86}, the Gonchar-Rakhmanov estimates in 
uniform rational approximation to 
sectionally holomorphic functions off an $S$-contour, and its 
generalization  to best $L^2$ and  $L^p$ approximants  in
\cite{BaYaSta,Totik}.

The present paper is, in part, a sequel to \cite{BaYaSta}. 
In the latter reference best $L^2$ and 
$L^\infty$ rational approximants are compared in the  $n$-th root sense,
whereas here we compare them in norm. 
We emphasize that the $L^2$ norm and
weighted variants thereof are of great importance in applications, due to 
their interpretation as a variance in a stochastic context. 
Moreover, best rational $H^2$ approximants have the 
interesting property of being attained through interpolation \cite{Levin}.
Note also that certain functions, like
Blaschke products, can be approximated in $H^2$-norm but not in 
the uniform norm by rational functions.

A key to the above-mentioned comparison is the derivation of lower bounds
on  the $L^2$ approximation error. Lower bounds in approximation 
are usually difficult to obtain; we dwell here on a topological machinery 
developed in \cite{BS02} which expresses the approximation error as the 
solution to a  $\max-\min$ problem, and we rely as well
on the Adamjan-Arov-Krein
theory of best uniform meromorphic approximation. We prove a somewhat general
result (Theorem~\ref{LBH2AAK}) which gives a lower bound on the $L^2$-best
rational approximation error of given degree, in terms of the ratios
of $L^2$ and $L^\infty$ norms of the singular vectors of the Hankel
operator with symbol the approximated function. We then apply it to three
cases where these ratios can be estimated: rational functions, Blaschke 
products, and Cauchy integrals on geodesic arcs. We use also
the $\max-\min$ principle to study linearized errors as a means to compute
further lower bounds. We also include numerical experiments, some of which
give excellent accuracy to estimate the $H^2$ error in rational 
approximation (see Table~\ref{results} in 
Section~\ref{sec:numericalresults}). To the author's knowledge, such results 
are first of their kind.

The paper is organized as follows. After some preliminaries on Hardy spaces
in Section~\ref{prelim}, we present in Section~\ref{Rat} the approximation
problems that we consider. Section~\ref{dual} is an  introduction
to the results
of \cite{BS02} and  it contains a basic  account of the  
Adamjan-Arov-Krein theory.
We derive in Section~\ref{Lbounds} our main theorem 
giving  lower bounds in $L^2$ rational 
approximation,
and we apply it to cases mentioned above. Finally, in section
\ref{LES}, we discuss linearized errors.

\section{Notations and Preliminaries}
\label{prelim}
Let $\DD$ be the unit disk in the complex plane $\CC$,  and 
$\TT$ the unit circle. We denote by $C(\TT)$ the space of continuous, 
complex-valued functions on $\TT$.
For $1\leq p\leq \infty$,
we put $L^p=L^p(\TT)$ for the familiar
Lebesgue space of complex measurable functions on $\TT$ 
such that
\[\|f\|_{p}=\left(\frac{1}{2\pi}
\int_{0}^{2\pi}|f(e^{i\theta})|^p\,d\theta\right)^{1/p}<\infty
\quad \text{if }1\leq p<\infty, \qquad
\|f\|_\infty={\rm ess.}\sup_{\theta \in 
[0,2\pi]}|f(e^{i\theta})|<\infty.\]

Hereafter, we let
$H^2=H^2(\DD)$ be the Hardy space of holomorphic functions in $\DD$ whose 
Taylor coefficients at 0 are square summable:
\[H^2=\{f(z)=\Sigma_{k=0}^\infty a_kz^k:\quad
\|f\|_{H^2}:=\Sigma_{k=0}^\infty |a_k|^2<+\infty\}.
\]
We refer the reader to \cite{gar} for standard facts on Hardy spaces.
By Parseval's relation 
\begin{equation}
\label{defintH2}\|f\|_{H^2}^2=\sup_{0\leq r<1}\frac{1}{2\pi}\int_0^{2\pi} 
|f(re^{i\theta})|^2\,d\theta,
\end{equation}
and  the map
\[\Bigl(f(z)=\Sigma_{k=0}^\infty a_kz^k\Bigr)\longrightarrow
\Bigl({f^*}(e^{i\theta}):=\Sigma_{k=0}^\infty a_ke^{ik\theta}\Bigr)
\]
is an isometry from $H^2$ onto the closed subspace of $L^2$
comprised of functions whose Fourier coefficients of strictly negative
index do  vanish. As is customary, we shall 
identify $H^2$ with this subspace so that the 
distinction  between $f$ and ${f^*}$ as well as  
$\|f\|_{H^2}$ and $\|f^*\|_2$ will disappear.
This conveniently allows one to regard members of
the Hardy class both as functions on $\DD$ and on $\TT$.
From the function-theoretic viewpoint, the correspondance $f\mapsto f^*$
is that ${f^*}(e^{i\theta})$ is almost everywhere the limit of 
$f(z)$ as $z$ tends non-tangentially to $e^{i\theta}$ within $\DD$.

We put $\bar{H}^{2,0}=\bar{H}^{2,0}(\CC\setminus\overline{\DD})$ 
for the companion Hardy space of holomorphic functions in 
$\CC\setminus\overline{\DD}$,  vanishing at infinity, whose 
Taylor coefficients there  are square summable:
\[\bar{H}^{2,0}=\{f(z)=\Sigma_{k=1}^\infty a_kz^{-k}:\quad
\|f\|_{\bar{H}^{2,0}}:=\Sigma_{k=1}^\infty |a_k|^2<+\infty\}.
\]
The map
\[\Bigl(f(z)=\Sigma_{k=1}^\infty a_kz^{-k}\Bigr)\longrightarrow
\Bigl({f^*}(e^{i\theta})=\Sigma_{k=1}^\infty a_ke^{-ik\theta}\Bigr)
\]
is an isometry from $\bar{H}^{2,0}$ onto the closed subspace of $L^2$
comprised of functions whose Fourier coefficients of non-negative index do 
vanish, and as before we identify $\bar{H}^{2,0}$ with the latter.
Clearly we have an orthogonal sum:
\begin{equation}
\label{decompP}
L^2=H^2\oplus\bar{H}^{2,0}.
\end{equation}

In fact, it holds that $f\in\bar{H}^{2,0}$ if and only if the function
$\check f$ given by
\begin{equation}
\label{defcheck}
\check{f}(z):=z^{-1}\overline{f(1/\bar{z})}
\end{equation}
 lies in $H^2$,
and the map  $f\mapsto \check{f}$ is an involutive
isometry of $L^2$ sending $H^2$ onto $\bar{H}^{2,0}$. Actually, $\check{f}$ has
same modulus as $f$ \emph{pointwise} on $\TT$ since  
$\overline{f(1/\bar{z})}=\overline{f(z)}$ when $|z|=1$. If $f$ 
is holomorphic on $\Omega$, then $f^\sharp(z)= 
\overline{f(1/\bar{z})}$ is holomorphic on the reflection of $\Omega$
across $\TT$, and if $f$ is rational $f^\sharp$ is likewise rational.
Of course, a relation like $f^\sharp=\bar{f}$
must be understood to hold on $\TT$ only.

We let
\[{\bf P}_+\bigl(\Sigma_{k\in\ZZ}\,\, a_ke^{ik\theta}\bigr)=
\Sigma_{k\geq0}\,\, a_ke^{ik\theta}\qquad\text{and}\qquad
{\bf P}_-\bigl(\Sigma_{k\in\ZZ} a_ke^{ik\theta}\bigr)=
\Sigma_{k<0} a_ke^{ik\theta}\]
indicate the so-called Riesz projections that discard the Fourier coefficients
of strictly negative and non-negative index respectively.
Clearly ${\bf P}_+$ (resp. ${\bf P}_-$)  contractively maps $L^2$
onto $H^2$ (resp.   $\bar{H}^{2,0}$) and ${\bf P}_++{\bf P}_-=I$.
We call ${\bf P}_+$ the \emph{analytic projection}
and ${\bf P}_-$ the \emph{anti-analytic projection}. Note that, 
by Cauchy's formula, ${\bf P}_\pm(f)$ can be expressed as Cauchy integrals:
\begin{equation}
\label{Cauchyproj}
{\bf P}_+(f)(z)=\frac{1}{2i\pi}\int_{\TT}\frac{f(\zeta)}{\zeta-z}d\zeta, \quad
|z|<1,\qquad
{\bf P}_-(f)(z)=\frac{1}{2i\pi}\int_{\TT}\frac{f(\zeta)}{z-\zeta}d\zeta, \quad
|z|>1.
\end{equation}

The Hardy space $H^\infty=H^\infty(\DD)$ consists of bounded holomorphic 
functions on $\DD$, endowed with the {\it sup} norm.
From \eqref{defintH2} we see that $H^\infty$ embeds contractively in $H^2$, 
in particular each $f\in H^\infty$ has
a non-tangential limit ${f^*}$ on $\TT$. It can be shown that
$\|{f^*}\|_\infty= \|f\|_{H^\infty}$, and that the map
$f\mapsto{f^*}$
is an isometry from $H^\infty$ onto the closed subspace of $L^\infty$
comprised of functions whose Fourier coefficients of strictly negative
index do  vanish. 
Again we identify $H^\infty$ with this subspace.
Likewise, the space $\bar{H}^{\infty,0}$ of bounded holomorphic functions 
vanishing at infinity in 
$\CC\setminus\overline{\DD}$ identifies {\it via} non-tangential limits
with the closed subspace of $L^\infty$
consisting of functions whose Fourier coefficients of non-negative
index do vanish. 
However, in contrast with the situation for $L^2$,
the operators ${\bf P}_{\pm}$ are unbounded on $L^\infty$.
Besides the norm topology, $H^\infty$ inherits the weak-* topology 
from $L^\infty(\TT)$. It is characterized by 
the fact that $f_n$ tends weak-* to $f$ if and only if
$\int_\TT f_n\varphi\to\int_\TT f\varphi$ for every $\varphi\in L^1$.
It is equivalent to require that $(\|f_n\|_\infty)_n$ is a bounded sequence and that,
for each $k$, the $k$-th Fourier coefficient of $f_n$ converges to the
$k$-th Fourier coefficient of $f$.

As is well-known \cite[ch. II, cor. 5.7]{gar},
a nonzero $f\in H^2$ factors uniquely
as $f=jw$ where
\begin{equation}
\label{defext}
w(z)=\exp\left\{\frac{1}{2\pi}\int_0^{2\pi}\frac{e ^{i\theta}+z}
{e^{i\theta}-z}\log|f(e ^{i\theta})|\, d\theta\right\}
\end{equation}
belongs to $H^2$ and is called the {\em outer factor} of $f$, normalized so as to be positive at zero,
while
$j\in H^\infty$ has modulus 1 a.e. on $\TT$ and is called the
{\em inner factor} of $f$. The latter may be further decomposed as
$j=bS$, where
\begin{equation}
\label{defBlaschke}
b(z)=cz^k\prod_{\zeta_l\neq0}\frac{{-\bar \zeta}_l}{|\zeta_l|}\,\frac{z-\zeta_l}{1-{\bar \zeta}_lz}
\end{equation}
is the normalized {\em Blaschke product}, with multiplicity $k\geq0$ at
the origin, associated to a sequence of points
$\zeta_l\in\DD\setminus\{0\}$
and to a constant $c\in\TT$, while
$$S(z)=\exp\left\{-\frac{1}{2\pi}\int_0^{2\pi}\frac{e ^{i\theta}+z}{e
    ^{i\theta}-z}\, d\mu(\theta)\right\}$$
is the {\em singular inner factor} associated with a positive 
singular measure $\mu$
on $\TT$. 
The $\zeta_l$ are of course the zeros of $f$ in $\DD$, counting
multiplicities by repetition. The number of zeros, finite or infinite, 
is called the degree of the Blaschke product. Throughout,
we let $B_n$ denote
the set of Blaschke products of degree at most $n$.
If the degree is infinite,
the convergence of the product in (\ref{defBlaschke})
is equivalent to the condition
\begin{equation}
\label{Bcond}
\sum_l(1-|\zeta_l|)<\infty
\end{equation}
 which holds automatically when $f\in H^2$. 
That $w(z)$ is well-defined rests on the fact that
$\log|f|\in L^1$ if $f\in H^2\setminus\{0\}$D.
A function $f\in H^2$ with inner-outer factorization $f=jw$
lies in $H^\infty$ if, and only if $w\in L^\infty(\TT)$.
For simplicity, we often say that a function is outer (resp. inner) 
if it is equal to its outer (resp. inner) factor.

We put $\cP_n[z]$ for the space of complex algebraic polynomials
of degree at most $n$ in the variable $z$, or simply $\cP_n$ if the
variable is understood. Below we let
$\mathcal{Z}(q)$ indicate the set of zeros of a polynomial $q$. 
For $q_n\in\cP_n[z]$, we define its
{\emph{reciprocal polynomial}} to be
\[\widetilde{q}_n(z):=z^{n}\,\overline{q_n(1/\bar z)}.\]
We warn the reader that this definition depends on $n$: if we consider
$q_{n-1}\in\cP_{n-1}$ as an element of $\cP_n$ with zero leading 
coefficient, the definitions of $\widetilde{q}_{n-1}(z)$ in $\cP_{n-1}$ and
in $\cP_{n}$ may be inconsistent. Therefore we always specify, {\it e.g.}
via a subscript ``$n$'' as in ``$q_n$'', which definition is used. Clearly the ``tilde''
operation is an involution of $\cP_n$ preserving modulus pointwise on $\TT$.

We designate by $\cR_{m,n}=\cR_{m,n}(z)$ 
the set of complex rational functions of type $(m,n)$ in
$L^2$, namely those that can be written as $p_m/q_n$
where $p_m$ belongs to $\cP_m$ and $q_n\in \cP_n$
has no root on ${\TT}$. 
When $r=p_m/q_n$  is in irreducible form,
the integer $\max\{m,n\}$ is the (exact) degree of~$r$.
Note that $B_m\subset\cR_{m,m}$ 
is comprised of rational functions of degree at most $m$ which are
analytic in $\DD$ and have unit modulus everywhere on $\TT$. 
Alternatively, $B_m$ consists of
functions $q_m/\widetilde{q}_m$  where 
$q_m\in \cP_m$ has all its roots in $\DD$.
Clearly, $B_m$ is included in the unit sphere of both $H^2$ and~$H^\infty$. 

We further set 
\[
H_m^2:=\{\frac{g}{q_m}:\ g\in H^2,\ q_m\in \cP_m\}.
\]
Members of $H_m^2$ identify in $L^2$ with non-tangential limits of 
meromorphic functions  with at most $m$ poles in $\DD$ 
(counting multiplicities) whose
$L^2$-means over $\{|z|=r\}$ remain eventually bounded 
as $r\to 1^-$. Functions in $\cup_m H_m^2$ are called
meromorphic in $L^2$.
Two equivalent descriptions of $H_m^2$ are useful:
on the one hand we get by pole-residue decomposition that
$H_m^2=H^2+(\cR_{m-1,m}\cap \bar H^{2,0})$, on the other hand we 
have that $H^2_m=B_m^{-1}H^2$, the set of quotients of $H^2$-functions by
Blaschke products of degree at most $m$.
Likewise we put
\[
H_m^\infty:=H^2_m\cap L^\infty=B_m^{-1}H^\infty=\{\frac{g}{q_m}:\ g\in H^\infty,\ q_m\in \cP_m\}
\]
for the set of meromorphic functions with at most $m$ poles in $L^\infty$.
\section{Best rational and meromorphic approximation in \texorpdfstring{$L^2$}{L2}}
\label{Rat}
For $n\geq1$ an integer, 
the best rational approximation problem of degree $n$ in $L^2$ is:

\noindent{\bf Problem R(n):} \emph{Given $h\in L^2$, to find $r^*\in \cR_{n,n}$ such that
\[\|h-r^*\|_2=\min_{r\in \cR_{n,n}}\|h-r\|_2.\]}

Write $h=h_1+h_2$ with $h_1\in H^2$, $h_2\in\bar{H}^{2,0}$.
By partial fraction expansion, each $r\in \cR_{n,n}$ can be decomposed as
$r_1+r_2$ where $r_1\in H^2$, $r_2\in\bar{H}^{2,0}$, and
$\mbox{deg\ }r_1+\mbox{deg\ }r_2\leq n$. Then, by (\ref{decompP}),
\[\|h-r\|_2^2=\|h_1-r_1\|_2^2+\|h_2-r_2\|_2^2\]
so that problem R($n$) reduces, modulo optimal allocation of the degrees
of $r_1$ and $r_2$ ($n+1$ choices), to a pair of problems of the
following types:

\noindent{\bf Problem RA(n):} \emph{Given $f\in H^2$, to find $r^*\in \cR_{n,n}\cap H^2$ 
such that
\[\|f-r^*\|_2=\min_{r\in \cR_{n,n}\cap H^2}\|f-r\|_2.\]}

\noindent{\bf Problem RAB(n):}
\emph{Given $f\in \bar{H}^{2,0}$, to find $r^*\in \cR_{n-1,n}\cap \bar{H}^{2,0}$ 
such that
\[\|f-r^*\|_2=\min_{r\in \cR_{n-1,n}\cap \bar{H}^{2,0}}\|f-r\|_2.\]}

In ``RA($n$)'' and ``RAB($n$)'',
the letter $"A"$ is mnemonic for ``analytic''
and  $"B"$ 
stands for ``bar''.

Problem RA($n$) is in fact equivalent to RAB($n$).
For we can parametrize  $r\in \cR_{n,n}\cap H^2$
as $r(0)+zr_3$ where $r(0)\in \RR$ and $r_3\in \cR_{n-1,n}\cap H^2$
vary independently, and by Parseval's theorem
\[\|f-r\|_2^2=|f(0)-r(0)|^2+\|(f-f(0))-zr_3\|_2^2\] 
hence $r(0)=f(0)$ is the optimal choice.
Thus, since multiplication by $1/z$ is an isometry,
we find upon replacing $f$ by $(f-f(0))/z$ 
that Problem RA($n$) is equivalent to the normalized version:

\noindent{\bf Problem RAN(n):} \emph{Given $f\in H^2$, to find 
$r^*\in \cR_{n-1,n}\cap H^2$ such that
\[\|f-r^*\|_2=\min_{r\in \cR_{n-1,n}\cap  H^2}\|f-r\|_2.\]}

Now, applying the check operation defined in \eqref{defcheck}, which
preserves $\cR_{n-1,n}$ and the degree, this last problem is 
seen to be equivalent to RAB($n$), as announced.
Note that when passing from RA($n$) to RAB($n$), the initial 
$f\in H^2$ to be approximated from $\cR_{n,n}\cap H^2$
gets transformed into  the function
$\overline{f(1/\bar z)}-\overline{f(0)}\in \bar{H}^{2,0}$ to be approximated from
$\cR_{n-1,n}\cap\bar{H}^{2,0}$.
Finally, we state the best meromorphic approximation problem 
with at most $n$ poles in $L^2$: 

\noindent{\bf Problem MA(n):} 
\emph{Given $f\in L^2$, to find $g^*\in H^2_n$ 
such that
\[\|f-g^*\|_2=\min_{g\in H^2_n}\|f-g\|_2.\]} 

Problem MA($n$) is also equivalent to
RAB($n$). Indeed, $H_n^2=H^2+(\cR_{n-1,n}\cap\bar{H}^{2,0})$
so that, by  orthogonality of 
$H^2$ and  $\bar{H}^{2,0}$, the $H^2$-component of a minimizer
in MA($n$) must be ${\bf P}_+(f)$ while the 
$\bar{H}^{2,0}$-component of this minimizer is a solution to
RAB($n$) with $f$ replaced by ${\bf P}_-(f)$. 

Let us mention that best meromorphic approximation,
unlike best rational approximation, is
conformally invariant. This makes it of independent interest in
a broader context, see \cite[prop. 5.4]{BMSW06}
for further details.
.

Having reduced all previous approximation problems to RAB($n$), hereafter
we discuss the latter.
It is known that RAB($n$) has a solution  which needs not be 
unique, and every solution has exact degree $n$ unless $f$ is rational 
of degree at most $n-1$ \cite{Erohin,Levin,Bar86}.

We shall write $d_2(f,\cR_{n-1,n})$ (resp. $d_2(f,\cR_{n,n})$)
for the distance from $f$ to $\cR_{n-1,n}$ (resp. $\cR_{n,n}$) in
$L^2$. For instance if  $f\in \bar{H}^{2,0}$, then $d_2(f,\cR_{n-1,n})$ is
both the value of  Problem RAB($n$) and of Problem MA($n$); and if
$f\in H^2$, then  $d_2(f,\cR_{n-1,n})$  (resp.  $d_2(f,\cR_{n,n})$)
is the value of problem RAN($n$) (resp. RA($n$)).
Besides, the value of MA($n$)  is denoted by $d_2(f, H^2_n)$.

When $f\in L^\infty$,
we let $d_{\infty}(f, H^\infty_n)$ indicate the distance from $f$ to 
$H^\infty_n$.
This is the value of the best meromorphic approximation problem with at most 
$n$ poles in $L^\infty$, that we did not formally introduce but which
stands analog
to MA($n$) with $L^2$ replaced by $L^\infty$ and $H^2_n$ by $H^\infty_n$.
We put also $d_\infty(f,\cR_{n-1,n})$ (resp. $d_\infty(f,\cR_{n,n})$)
for the distance from $f$ to $\cR_{n-1,n}$ (resp. $\cR_{n,n}$)  in $L^\infty$.
\section{Duality in meromorphic approximation}
\label{dual}
Pick $f\in \bar{H}^{2,0}$ and let us parametrize 
$r\in \cR_{n-1,n}\cap \bar{H}^{2,0}$ as $r=p_{n-1}/q_n$
where $p_{n-1}$ ranges over $P_{n-1}$ and $q_n$ ranges over those polynomials
in $\cP_n$ whose roots lie in $\DD$.
Then $q_n/\widetilde{q}_n\in B_n$ and since 
$p_{n-1}/\widetilde{q}_n\in H^2$ we have by
orthogonality of  $H^2$ and  $\bar{H}^{2,0}$ that
\begin{equation}
\label{multB}
\|f-\frac{p_{n-1}}{q_n}\|_2^2=\|f\frac{q_n}{\widetilde{q}_n}-
\frac{p_{n-1}}{\widetilde{q}_n}\|_2^2=
\|{\bf P}_-(f\frac{q_n}{\widetilde{q}_n})\|_2^2+
\|{\bf P}_+(f\frac{q_n}{\widetilde{q}_n})
-\frac{p_{n-1}}{\widetilde{q}_n}\|_2^2.
\end{equation}
Clearly the product of a $\bar{H}^{2,0}$-function by a polynomial in $\cP_n$
yields a member of $z^n\bar{H}^{2,0}$. Therefore
\begin{equation}
\label{minp}
\widetilde{q}_n{\bf P}_+(f\frac{q_n}{\widetilde{q}_n})=
fq_n-\widetilde{q}_n{\bf P}_-(f\frac{q_n}{\widetilde{q}_n})\in 
z^n\bar{H}^{2,0}\cap H^2=\cP_{n-1},
\end{equation}
entailing that $p_{n-1}=\widetilde{q}_n{\bf P}_+(fq_n/\widetilde{q}_n)$
is the minimizing choice in \eqref{multB} for fixed $q_n$.
Consequently 
\begin{equation}
\label{vpreHank}
\min_{r\in \cR_{n-1,n}\cap \bar{H}^{2,0}}\|f-r\|_2=
\min_{q_n\in \cP_n, \mathcal{Z}(q_n)\subset\DD}
\|{\bf P}_-(f\frac{q_n}{\widetilde{q}_n})\|_2=
\min_{b_n\in B_n}
\|{\bf P}_-(fb_n)\|_2.
\end{equation}
That the infimum is indeed attained in the right hand side of \eqref{vpreHank}
follows from \eqref{multB} and the fact that RAB($n$) has a solution.
Define $A_f$, the {\em Hankel operator with symbol $f$},
by 
\begin{equation}
\label{Hankel2}
\begin{array}{lll}
A_{f}:  H^{\infty} & \longrightarrow &  \bar{H}^{2,0} \\
\ \ \ \ \ \ \      v &          \mapsto & {\bf P}_-(fv). \\
\end{array}
\end{equation}
It is evident that $A_f$ is continuous and that $|||A_f|||=\|f\|_2$,
a unit maximizing vector being $v\equiv 1$.
Here and below, we let $|||.|||$ stand for the 
operator norm, and a maximizing vector of an operator $E$ is a nonzero vector
$v$ such that $\|E v\|/\|v\|=|||E|||$. 

The content of the discussion leading from \eqref{multB} to  \eqref{vpreHank} 
may now be restated as follows.
\begin{prop}
\label{lemH}
For $f\in \bar{H}^{2,0}$, it holds that
\begin{equation}
\label{minHankB}
d_2(f,\cR_{n-1,n})=\min_{b_n\in B_n}\|A_f(b_n)\|_2.
\end{equation}
A rational function $p_{n-1}/q_n\in\cR_{n-1,n}$ is a solution to 
\textrm{RAB}($n$) if, and only if $b_n=q_n/\widetilde{q}_n$ 
is a minimizing Blaschke product in
\eqref{minHankB} and
$p_{n-1}=\widetilde{q}_n{\bf P}_+(fb_n)$.
\end{prop}

Put $\mathcal{L}_k$ for the space of linear operators
from $H^\infty$ into $\bar{H}^{2,0}$ which are weak-* continuous
and have rank not exceeding $k$. 
For $k=0,1,2,...$, we denote by $\sigma_k(A_f)$  
the $k$-th {\em  approximation
number} of $A_f$ defined by
\begin{equation}
\label{singdef}
\sigma_k(A_f)=\inf \bigl\{ 
|||A_f-\Gamma|||,~~\Gamma\in\mathcal{L}_k.\bigr\}.
\end{equation}
Note that $\sigma_k(A_f)\geq \sigma_{k+1}(A_f)$  and that 
$\sigma_0(A_f)=|||A_f|||$.

We need also introduce the {\emph{genus}}
of a closed symmetric subset $K$
in a topological vector space; here, 
\emph{symmetric} means that if $v\in K$ then also $-v\in K$.
By definition the genus of $K$,  denoted by ${\bf gen} (K)$,
is the smallest positive integer $m$
for which there exists an odd continuous mapping 
\begin{equation}
\label{defgen}
G: K \longrightarrow \RR^m\setminus\{0\},
\end{equation}
or else $+\infty$ if no finite $m$ meets the above requirement.
By convention the genus is zero if $K=\emptyset$. When $K$ is 
compact  and does not contain $0$, then
${\bf gen} (K)$ is always finite, see \cite{Zeidler}.
For instance, if $m\geq1$, 
the classical Borsuk-Ulam theorem from topology 
\cite[ch. 2, sec. 6]{GuilleminPollack} implies
that any symmetric set in $\RR^m$ which is homeomorphic to the 
(real) $(m-1)$-dimensional
Euclidean sphere $\sph^{m-1}$ through an odd map has genus $m$. 

Below, we shall be concerned  with weak-* compact subsets of
$\mathcal{S}^\infty$, the unit sphere of $H^\infty$. In this connection,
we let
$$ {\cal K}_m^\infty=\bigl\{K \subset{\cal S}^\infty : K \text{ is a weak-* compact 
symmetric subset of }{\cal S}^\infty\text{ with {\bf gen}}(K) \geq m \bigr\}.$$

Subsequently, we define the (generalized) \emph{singular numbers} of $A_f$
by
\begin{equation}
\label{minmax}
\lambda_m(A_f)
=\max_{K \in {\cal K}_m^\infty} \min_{u \in K}\|A_f(u)\|_2,\qquad m=0,1,2,...
\end{equation}
The following theorem, which was established in \cite{BS02},
connects approximation numbers and singular numbers of $A_f$ with
the value of Problem RAB($n$):
\begin{theorem}{\rm \cite[thm. 8.1]{BS02}}
\label{singlust}
Let $f\in \bar{H}^{2,0}$ and $A_f:H^\infty\to \bar{H}^{2,0}$ 
the Hankel operator with symbol $f$. 
For each integer $n\geq0$, the following equalities hold:
\begin{equation}
\label{critsinglust}
d_2(f,\cR_{n-1,n})=\sigma_n(A_f)=\lambda_{2n+1}(A_f)=\lambda_{2n+2}(A_f).
\end{equation}
\end{theorem}

Theorem~\ref{singlust} is reminiscent of  a famous theorem by
Adamjan-Arov-Krein  
(in short: the AAK theorem) characterizing $d_{\infty}(f,H^\infty_n)$ 
rather than
$d_2(f,\cR_{n-1,n})$. To state the result, let us define for $f\in L^\infty$ the Hankel operator 
$\Gamma_f$ by 
\begin{equation}
\label{Hankelinf}
\begin{array}{lll}
\Gamma_{f}:  H^{2} & \longrightarrow &  \bar{H}^{2,0} \\
\ \ \ \ \ \ \      v &          \mapsto & {\bf P}_-(fv). \\
\end{array}
\end{equation}
Although the definitions of $A_f$ and $\Gamma_f$ are formally the same,
observe that the domains in \eqref{Hankel2} and \eqref{Hankelinf}
are different.
The definition of $s_k(\Gamma_f)$ is still given by \eqref{singdef} except that
$A_f$ is replaced by $\Gamma_f$ and $\Gamma$ now ranges over
linear operators from $H^2$ into $\bar{H}^{2,0}$ having 
rank at most $k$. If in addition $f$ is continuous on $\TT$, 
then $\Gamma_f$ is
compact \cite[ch. 1, thm. 5.5]{Peller}. Then, if we let $\Gamma_f^*$ 
denote the adjoint,
$\Gamma_f^*\Gamma_f$ is a compact selfadjoint
operator from the Hilbert space $H^2$ into itself and as such 
it has a complete orthonormal family of eigenvectors called the singular 
vectors of $\Gamma_f$; the associated eigenvalues are none 
but the squared approximation numbers of $\Gamma_f$ 
\cite[ch. II, thm. 2.1]{GK}, and there holds 
the Courant $\max\min$ principle \cite[sec. 22.11a]{ZeidlerII}:
\begin{equation}
\label{CPHilb}
s_n(\Gamma_f)=\max_{V\in\mathcal{V}_{n+1}} \min_{\stackrel{v \in V}{\|v\|_2=1}}\|\Gamma_f(v)\|_2,
\end{equation}
where $\mathcal{V}_{n+1}$ is the collection of linear subspaces of $H^2$
of complex dimension at least $n+1$.
In this Hilbertian context, the approximation number $s_n(\Gamma_f)$ is also called the $n$-th singular value of $\Gamma_f$. We say that a function $v$ is \emph{associated with a singular value $s$} when $v$ is an eigenvector of $\Gamma_f^*\Gamma_f$ associated with the eigenvalue $s^2$: $v = s^2\Gamma_f^*\Gamma_f(v)$. As a particular case of Equation~\eqref{CPHilb} a maximizing vector is
just a singular vector associated with~$s_0(\Gamma_f)$.
\begin{theorem}[The AAK theorem]{\rm \cite[thms. 0.1 \& 0.2]{AAK}\cite[ch. 4, thm. 1.2]{Peller}}
\label{AAKth}
Let $f\in L^\infty$ and $\Gamma_f:H^2\to \bar{H}^{2,0}$ 
be the Hankel operator with symbol $f$. 
For each integer $n\geq0$, it holds that
\begin{equation}
\label{egAAK}
d_{\infty}(f,H^\infty_n)=s_n(\Gamma_f).
\end{equation}
If in addition
$f\in C(\TT)$, then $\Gamma_f$ is compact and the quantity
\eqref{egAAK} is also equal to \eqref{CPHilb}.
\end{theorem}
The case $n=0$ of Theorem~\ref{AAKth}, {\it i.e.} 
that  $|||\Gamma_f|||=d_{\infty}(f,H^\infty)$ was known earlier as Nehari's theorem.

If we compare 
\eqref{minmax} and \eqref{critsinglust}
with \eqref{CPHilb} and \eqref{egAAK}
for $f\in \bar{H}^{2,0}\cap L^\infty$, we see that
the main difference 
between best meromorphic approximation with at most $n$ poles in
$L^2$ and in $L^\infty$  lies with the
maximization step in \eqref{minmax}, which in the $L^2$-case  must
be taken over all compact sets of genus at least\footnote{
That $\lambda_{2n+1}(A_f)=\lambda_{2n+2}(A_f)$ in \eqref{critsinglust} is 
inessential and due the fact that $A_f$ is complex linear whereas the genus 
is a real notion.} $2n+2$ and not just
Euclidean spheres of real dimension $2n+1$. 
It follows from  \cite[thm 1]{Bar90} or \cite[thm. 5.3]{HP} 
that $B_n$ is homeomorphic to $\sph^{2n+1}$
and inspection of the proof reveals that the homeomorphism is odd.
Moreover  $B_n$ is weak-* compact
in $\mathcal{S}^\infty$ \cite[lem. 7.3]{BS02},
therefore $B_n\in\mathcal{K}_{2n+2}^\infty$
and from Proposition~\ref{lemH} we see that it is
a supremizer in \eqref{minmax}.

We mention for completeness
a companion to Theorem~\ref{singlust} dealing with
$\min\max$ (not $\max\min$):
\begin{theorem}{\rm \cite{Prokhorov2002}\cite[eqn. (78)]{BS02}}
\label{singvect}
Let $f\in \bar{H}^{2,0}$ and $A_f:H^\infty\to \bar{H}^{2,0}$ 
be the Hankel operator with symbol $f$. 
For each integer $n\geq0$, the following equality hold:
\begin{equation}
\label{critsinglustVa}
d_2(f,\cR_{n-1,n})= \min_{W\in{\cal W}_n}
\max_{\stackrel{w\in W}{\|w\|_{\infty}= 1}} 
\,  \|A_f(v)\|_2,
\end{equation}
where $\mathcal{W}_n$ is the collection of linear subspaces 
in $H^\infty$ of (complex) codimension at most $n$. 
\end{theorem}
Note that \eqref{critsinglustVa} is the exact counterpart for $A_f$ of
the standard Courant $\min\max$ principle for $\Gamma_f$:
\[
d_{\infty}(f,H^\infty_n)= \min_{X\in{\cal X}_n}
\max_{\stackrel{w\in X}{\|w\|_{2}= 1}} 
\,  \|\Gamma_f(v)\|_2,
\]
where   $\mathcal{X}_n$ is the collection of linear subspaces 
in $H^2$ of (complex) codimension at most $n$.

Using  Proposition~\ref{lemH} it is easy to see that if $p_{n-1}/q_n$ 
is a solution to RAB($n$),
then the subspace $(q_n/\widetilde{q}_n) H^\infty$, comprised of multiples of
$q_n/\widetilde{q}_n$ in $H^\infty$, is
a minimizing $W$  in \eqref{critsinglustVa}. 
In the rest of the paper, we use the maximizing step in \eqref{minmax} 
together with Theorem~\ref{singlust} to derive lower bounds for 
Problems RAB($n$).

\section{Lower bounds}
\label{Lbounds}
\subsection{Comparing \texorpdfstring{$L^2$}{L2} and \texorpdfstring{$L^\infty$}{L-infinity} meromorphic approximation}
Consider $f\in \bar{H}^{2,0}\cap L^\infty$ and $r,r^*\in\cR_{n-1,n}$ with
$r^*$  a solution to RAB($n$), {\it i.e.} a best approximant to $f$ in $L^2$
from $\cR_{n-1,n}$. Then $\|f-r^*\|_2\leq\|f-r\|_2$. Now, for any $h\in H^\infty$, Parseval's theorem gives $\|f-r\|_2 \le \|f-r-h\|_2$.
Finally, since the $L^\infty$-norm
dominates the $L^2$-norm $\|f-r-h\|_2\leq\|f-r-h\|_\infty$ and so we have
 \[\|f-r^*\|_2 \leq \|f-(r+h)\|_\infty.
\]
Thus, minimizing  over $r,h$, we find that
$d_2(f,\cR_{n-1,n})\leq d_{\infty}(f,H^\infty_n)$.
However, it is {\it a priori} unclear how large the gap between the two 
errors can be. Below,
dwelling on Theorems~\ref{singlust} and~\ref{AAKth}, we derive  
when $f$ is continuous a lower bound in terms of  the ratio between
$L^2$ and $L^\infty$ norms of the singular vectors of the Hankel operator 
$\Gamma_f$. 
\begin{theorem}
\label{LBH2AAK}
Let $f\in \bar{H}^{2,0}\cap C(\TT)$ and $n\geq0$ an integer. Consider
an orthonormal family $v_0,\cdots,v_n$ of singular vectors  of the
Hankel operator $\Gamma_f$ ({\it cf.} \eqref{Hankelinf}),
where $v_k$ is associated to the singular value $s_k(\Gamma_f)$.
Define $M_n(f):=\min\{d_{\infty}(f,H_j^\infty)/\|v_j\|_{\infty},\,0\leq j\leq n\}$
if $v_j\in H^\infty$ for $0\leq j\leq n$, and $M_n(f):=0$ otherwise.
Then
\begin{equation}
\label{borninfgen}
\frac{M_n(f)}{\sqrt{n+1}}\leq d_2(f,\cR_{n-1,n}).
\end{equation}
\end{theorem}
{\it Proof:} 
if $M_n(f)=0$, then \eqref{borninfgen} is trivial.
Otherwise, the linear span of  $\{v_0,\cdots,v_n\}$ over $\CC$ is a real 
$2n+2$-dimensional
vector space in $L^2\cap L^\infty$, and we may endow it either  with the 
$L^2$-norm or else with the
$L^\infty$-norm. Let $S_2$ and $S_\infty$ indicate the corresponding unit 
spheres. Identifying a vector with its coordinates,
we see that $S_2$ is just
$\sph^{2n+1}$, and clearly $v\mapsto v/\|v\|_\infty$ is an odd homeomorphism 
from $S_2$ onto $S_\infty$. Therefore, by the Borsuk-Ulam theorem,
$S_\infty$ is a compact set  of genus $2n+2$. Now, if we let
$v\in S_\infty$ and write $v=\sum_{j=0}^n \lambda_jv_j$ while abbreviating
$s_j(\Gamma_f)$ as $s_j$, we get using ``$\langle\,,\,\rangle$'' to mean
Hermitian scalar product on $\TT$ that
\begin{align}
\nonumber
\|A_f(v)\|_2^2&=\langle A_f (v),A_f (v)\rangle
=\langle \Gamma_f (v),\Gamma_f (v)\rangle=
\langle \Gamma^*_f\Gamma_f v,v\rangle=\Sigma_{j=0}^n|\lambda_j|^2s_j^2\\
\label{elcalc}
&\geq\frac{1}{n+1}\left(\sum_{j=0}^n |\lambda_j|s_j\right)^2
\geq\frac{M_n^2(f)}{n+1}\left(\sum_{j=0}^n|\lambda_j|\|v_j\|_\infty\right)^2
\geq \frac{M_n^2(f)}{n+1},
\end{align}
where the second line in \eqref{elcalc} uses the Schwarz inequality, the definition of $M_n(f)$ together with the equality $s_j(\Gamma_f)=d_{\infty}(f,H_j^\infty)$
from Theorem~\ref{AAKth}, the triangle inequality and the fact 
that $\|v\|_\infty=1$.
Inequality \eqref{borninfgen} now follows from \eqref{elcalc} and 
Theorem~\ref{singlust}.

\hfill\boite 

The kernels $\text{Ker} A_f$ and $\text{Ker} \Gamma_f$ are closed subsets of
$H^\infty$ and $H^2$ respectively, and clearly 
$\text{Ker} A_f=\text{Ker} \Gamma_f\cap H^\infty$.
({\it cf.} definitions \eqref{Hankel2} and \eqref{Hankelinf}).
By a theorem of Beurling \cite[ch. II, thm. 7.1]{gar},
being closed and shift-invariant ({\it i.e.} invariant under multiplication 
by the variable $z$), 
$\text{Ker} \Gamma_f$ is
either trivial ($\{0\}$ or $H^2$) or else 
consists of all multiples of some inner 
function $\mathfrak{j}$, that is, $\text{Ker} \Gamma_f=\mathfrak{j} H^2$. 
In the latter case 
$\text{Ker} A_f=\mathfrak{j} H^\infty$, in particular 
$\text{Ker} \Gamma_f$ and 
$\text{Ker} A_f$ are simultaneously nontrivial. In this
situation the proof of Theorem~\ref{LBH2AAK}
quickly leads to an improvement of itself as follows. 
Notations and 
assumptions being as in the theorem,
set $\|v_j\|_{H^\infty/\text{Ker}A_f}$ 
to be $+\infty$
if $v_j\notin H^\infty$  and to be the distance 
from $v_j$ to $\text{Ker}A_f$ in $H^\infty$ otherwise.  Observe that if
$\|v_{j_0}\|_{H^\infty/\text{Ker}A_f}=0$ for some $j_0\in\{0,\cdots,n\}$,
then $v_{j_0}\in \text{Ker} \Gamma_f$ 
which entails
that $\Gamma_f$ has rank 
at most $j_0$ by definition of singular values. It is
a theorem of Kronecker  \cite[ch. 1, cor. 3.2]{Peller}
that this happens if and only if 
$f\in H^\infty_{j_0}$, and since  $f\in \bar{H}^{2,0}\cap C(\TT)$ 
by assumption 
we get that $f\in \cR_{j_0-1,j_0}$. In particular it holds in this case that
$d_{\infty}(f,H^\infty_j)=d_2(f,\cR_{j-1,j})=\|v_{j}\|_{H^\infty/\text{Ker}A_f}=0$
for all $j\geq j_0$.  Keeping this observation in mind, let us define
\begin{equation}
\label{defQn}
Q_n(f):=\min_{0\leq j\leq n}\left\{\frac{d_{\infty}(f,H_j^\infty)}{\|v_j\|_{H^\infty/\text{Ker}A_f}}\right\}, 
\end{equation}
where $Q_n(f)$ is to be interpreted as $0$ 
if $\|v_{j_0}\|_{H^\infty/\text{Ker}A_f}=0$ for some $j_0\in\{1,\cdots,n\}$
(in which case $d_{\infty}(f,H_{j_0}^\infty)=0$ as well by what precedes).
\begin{corollary}
\label{corquot} 
Theorem~\ref{LBH2AAK} remains valid if $M_n(f)$ gets replaced by $Q_n(f)$.
 \end{corollary}
{\it Proof:} 
we can assume that $v_j\in H^\infty\setminus\text{Ker} A_f$ 
for $0\leq j\leq n$, otherwise 
$Q_n(f)=0$ and there is nothing to prove. By the discussion before the 
corollary, this amounts to say that $f\notin H^\infty_n$.
Next, pick $\varepsilon>0$ and $g_j\in\text{Ker} A_f$ such that
$\|v_j-g_j\|_\infty<\|v_j\|_{H^\infty/\text{Ker}A_f}+\varepsilon$ for 
each $j\in\{1,\cdots,n\}$. If we let $w_j=v_j-g_j$, then
$A_f(w_j)=\Gamma_f(w_j)=\Gamma_f(v_j)$ and the $w_j$ are linearly 
independent over $\CC$. Indeed, if $\sum_{j=0}^n\lambda_j w_j=0$ with
$\lambda_{j_0}\neq0$, applying $\Gamma^*_f\Gamma_f$ yields  
$\sum_{j=0}^n\lambda_j s^2_j(\Gamma_f) v_j=0$ and since the $v_j$ are 
linearly independent we
have that $s_{j_0}(\Gamma_f)=0$; thus, by the AAK theorem, we get 
that $f\in H^\infty_{j_0}\subset H^\infty_n$, 
contrary to our initial assumption.
Replacing now $v_j$ by $w_j$ in the proof of
Theorem~\ref{LBH2AAK} and using that $\Gamma_f(w_j)=\Gamma_f(v_j)$,  
we obtain instead of \eqref{elcalc} that, whenever
$w=\sum_{j=0}^n\lambda_jw_j$ is such that $\|w\|_\infty=1$, then
\[
\|A_f(w)\|_2^2\geq\frac{1}{n+1}\,\,\min_{0\leq j\leq n}\left(
\frac{d_{\infty}(f, H^\infty_j)}{\|v_j\|_{H^\infty/\text{Ker} A_f}+\varepsilon}
\right)^2.
\]
Thus, letting $\varepsilon$ go to $0$, we get the desired result from
Theorem~\ref{singlust} again.

\hfill\boite

Theorem~\ref{LBH2AAK} is useful only if we have a fair appraisal of $M_n(f)$.
The latter is delicate to estimate in general, but in the following
subsections we point out three cases where this can be done in different
guises. They are: the case of a general
rational function which can be approached numerically; the case of a Blaschke 
product where estimates  can be given in terms of the zeros;
the case of Cauchy integrals
over hyperbolic geodesic arcs
in which boundedness of $M_n(f)$ can be proved {\it via}  a careful analysis
of formulas behind AAK theory,  dwelling on the work in
\cite{BS02}.
\subsection{Application to rational functions}
\label{sec:applicationRational}
When $f$ is rational, the bounds in 
Corollary~\ref{corquot} can be numerically computed.
As explained in Section~\ref{Rat}, the general case reduces  by
partial fraction extension to the special case where $f\in \bar{H}^{2,0}$,
the detail of which is carried out below.

Write $f=p/q$ where $p\in \cP_{N-1}$, 
$q\in\cP_N$ is monic with all roots in $\DD$, and $p$, $q$ are 
coprime as polynomials. Let us write
\[q(z)=\Pi_{k=1}^N (z-\zeta_k)\]
where each $\zeta_k\in\DD$ is repeated according to multiplicity.
It is clear from definition \eqref{Hankelinf} 
that $\text{Ker}\Gamma_f$ consists of those $H^2$-functions
vanishing at the zeros of $q$, hence
$\text{Ker}\Gamma_f=(q/\widetilde{q})H^2$. Its orthogonal complement in
$H^2$ is $(\text{Ker}\,\Gamma_{f})^\perp=\cP_{N-1}/\widetilde{q}$, 
an orthonormal basis of which is given 
according to the Malmquist -Walsh lemma  by the formulas
\cite[ch. V, sec 1]{Nikolskii}:
\begin{equation}
\label{defMW}
e_j(z)=\frac{\left(1-|\zeta_j|^2\right)^{1/2}}{1-\bar{\zeta}_j z}
\,\,\Pi_{k=0}^{j-1}\,\frac{z-\zeta_k}{1-\bar{\zeta}_kz},\qquad
1\leq j\leq N,
\end{equation}
where the empty product is understood to be 1. The effect of $\Gamma_f$ on 
any member of  $(\text{Ker}\,\Gamma_{f})^\perp$ is  easily computed upon 
introducing $a\in\cP_{N-1}$ and $b\in\cP_{N-1}$ such that the following 
Bezout relation holds: $a\widetilde{q}+bq=1$. Indeed, one has for any 
$u\in {\cP}_{N-1}[z]$ that
\begin{equation}
\label{calcHankrat}
\Gamma_f(u/\widetilde{q})={\bf P}_-\left(\frac{pu}{q\widetilde{q}}\right)
={\bf P}_-\left(\frac{pua}{q}+\frac{pub}{\widetilde{q}}\right)
={\bf P}_-\left(\frac{pua}{q}\right)=\frac{R_q(pua)}{q},
\end{equation} 
where we used that $pub/\widetilde{q}\in H^2$ and, for any polynomial $P$, 
$R_q(P)$ indicates the remainder of Euclidean 
division of $P$ by $q$. In particular, we get from
\eqref{calcHankrat} that $\text{Im}\Gamma_f=\cP_{N-1}/q$. 
The Hermitian scalar product on $\TT$
can be computed in several ways for functions in  $\cP_{N-1}/q$;
one which does not use partial fraction expansion is as follows. 
Pick $u,v\in\cP_{N-1}$. 
Observing that $z^N/\widetilde{q}$ is conjugate to $1/q$ on $\TT$ and 
denoting with
$Q_q(P)$ the quotient of Euclidean division of the polynomial $P$ by $q$
(so that $P=qQ_q(P)+R_q(P)$), we get since  $a\widetilde{q}+bq=1$ that
\begin{equation}
\label{calcsp}
\langle \frac{u}{q}\,,\,\frac{v}{q}\rangle=
\langle \frac{z^Nu}{\widetilde{q}q}\,,\,v\rangle
=\langle\frac{z^Nub}{\widetilde{q}}+Q_q\left(z^Nua\right)+
\frac{R_q(z^Nua)}{q}\,,\,v\rangle
=\langle Q_q\left(z^Nua\right)\,,\,v\rangle,
\end{equation}
where we used that $R_q(z^Nua)/q\in\bar{H}^{2,0}$ and 
$z^Nub/\widetilde{q}\in z^N H^2$  are
both orthogonal to $v\in\cP_{N-1}$ by Parseval's theorem.
The last term in \eqref{calcsp} is now a scalar product between 
polynomials which can be computed as a Euclidean one in the basis 
$\{z^k;\ 0\leq k\leq N-1\}$.

Writing $e_j=u_j/\widetilde{q}$ where
$e_j$ was defined in \eqref{defMW}, we can use \eqref{calcHankrat},
\eqref{calcsp} to compute the Hermitian matrix
$M=\langle \Gamma^*_f\Gamma_f(e_i)\,,\,e_j\rangle=
\langle \Gamma_f(e_i),\,\Gamma_f(e_j)\rangle$, and
an orthonormal family of singular vectors $v_0,\cdots,v_{N-1}$ associated with the nonzero 
singular values of $\Gamma_f$ is then obtained by diagonalization of $M$
(of course any other orthonormal basis of
$\cP_{N-1}/\widetilde{q}$ than $(e_k)$ could be used as well). 
More precisely, the $k$-th 
row of a unitary matrix $U$ such that $UMU^*$ is diagonal yields
coordinates for $v_k$ in the basis $e_j$.
The diagonal terms are the squared singular values $s^2_k(\Gamma_f)$ for 
$0\leq k\leq N-1$, which are none but the $d_{\infty}(f,H^\infty_k)$ by the 
AAK theorem.
Moreover, it follows from Nehari's theorem that
\begin{equation}
\label{ecrasev}
\|v_j\|_{H^\infty/\text{Ker}A_f}=d_{\infty}(v_j\widetilde{q}/q,H^\infty)=
|||\Gamma_{v_j\widetilde{q}/q}|||,
\end{equation}
and the last term in \eqref{ecrasev} is the largest singular value of a 
Hankel operator with rational symbol which can be computed in the same 
manner as indicated above to compute  $s_0(\Gamma_f)$. 

Thus, we can evaluate $Q_n$ defined in  \eqref{defQn} for all $n$, 
hence also the lower bound on $d_2(f,\cR_{n-1,n})$ given by Theorem~\ref{LBH2AAK}
and Corollary~\ref{corquot}.
We implemented a prototype algorithm to compute these two bounds. Numerical experiments are presented in Section~\ref{sec:numericalresults}.

\subsection{Application to Blaschke products}
In this section, we use Theorem~\ref{LBH2AAK} to
derive some lower bounds for Problem RA($n$) when
$f$ is a Blaschke product of finite or infinite degree. This last case is 
instructive to contrast rational approximation in $L^2$ and $L^\infty$ norms, 
for on the one hand 
the value of Problem RA($n$) tends to
zero as $n$ goes large (since rational functions are dense in $H^2$), while
on the other hand  $f$ cannot be approximated ``at all'' 
by rational functions in $H^\infty$, {\it i.e.} zero is a best uniform
approximant. This follows from the lemma below which is not easy
to locate in the literature.
\begin{lem}
\label{bab}
Let $b$ be a Blaschke product and $n$ be a positive integer which is strictly 
less 
than the degree of $b$ (if $b$ has infinite degree the assumption is void).
Then
\begin{equation}
\label{impRA}
d_\infty(b,\cR_{n,n})=\|b\|_\infty=1.
\end{equation}
\end{lem}
{\it Proof:}  clearly $d_\infty(b,\cR_{n,n})\leq1$ for zero is a candidate 
approximant.
Moreover, if $r\in \cR_{n,n}\cap H^{\infty}$ then
$\bar{r}\in H^\infty_n$. Therefore,  upon conjugating, we get
$d_\infty(b,\cR_{n,n})\geq d_{\infty}(\bar{b},H^\infty_n)$
and it is enough to show the latter is at least 1, hence in fact
equal to 1.

Assume first that $b$ has finite degree $d$, and write
$b=q_d/\widetilde{q}_d$ where $q_d\in\cP_d$ has zeros in $\DD$ only.
Then $\bar{b}=\widetilde{q}_d/q_d$, and the kernel of $\Gamma_{\bar{b}}$ is
$b H^2$ whose orthogonal complement in $H^2$ is 
$(\text{Ker}\,\Gamma_{\bar{b}})^\perp=\cP_{d-1}/\widetilde{q}_d$
as pointed out in the previous section.
Now, if
$p_{d-1}\in\cP_{d-1}$, then $\Gamma_{\bar{b}}(p_{d-1}/\widetilde{q}_d)=
p_{d-1}/q_d$ so that $\Gamma_{\bar{b}}$ is an isometry from 
$(\text{Ker}\,\Gamma_{\bar{b}})^\perp$ onto its image.
Consequently the first $d$ singular values of $\Gamma_{\bar{b}}$ 
are equal to 1 (the remaining ones being zero). That 
 $d_{\infty}(\bar{b},H^\infty_n)=1$ now follows from the AAK theorem and the fact that
$n\leq d-1$.

Assume next that $b$ has infinite degree. We can write $b=b_{n+1}b_\infty$
where $b_{n+1}$ has degree $n+1$ and $b_\infty$ has infinite degree.
If $g \in H^{\infty}_n$ then also $b_\infty g\in H^\infty_n$, 
and since $|b_\infty|=1$ a.e. on $\TT$, we get by the first part of the proof 
that
\begin{equation}
\label{eqt}
\|\bar{b}-g\|_\infty=\|b_\infty\bar{b}-b_\infty g\|_\infty=\|\bar{b}_{n+1}-b_\infty g\|_\infty\geq1,
\qquad g\in H^\infty_n,
\end{equation}
hence  $d_{\infty}(\bar{b},H^\infty_n)\geq1$, as desired.

\hfill\boite

We turn to the main result of this section:
\begin{theorem}
\label{LBB}
Let $b$ be a Blaschke product, of finite or infinite degree. Let us arrange
its zeros into a (finite or infinite) sequence 
$\zeta_1,\zeta_2,\cdots,$  where each $\zeta_j$ is
repeated according to its multiplicity and the corresponding sequence of 
moduli is nondecreasing:
$|\zeta_1|\leq|\zeta_2|\leq\cdots$. For each positive integer
$n$ strictly less
than the degree of $b$ (if $b$ has infinite degree the assumption is void),
it holds that
\begin{equation}
\label{borninfB}
\frac{\left(1-|\zeta_{n+1}|^2\right)^{1/2}}{\sqrt{n+1}}\leq 
d_2(b,\cR_{n,n})
\end{equation}
and also that 
\begin{equation}
\label{borninfBa}
\left(\sum_{j=0}^{n}\frac{1}{\left(1-|\zeta_{j}|^2\right)^{1/2}}\right)^{-1}
\leq 
d_2(b,\cR_{n,n}).
\end{equation}
\end{theorem}
{\it Proof:} assume first that $b$ has finite degree $d$, so that 
$b\in C(\TT)$, and write
$b=q_d/\widetilde{q}_d$ where $q_d\in\cP_d$ has zeros in $\DD$ only.
By the equivalence between Problem RA($n$) and RAN($n$) 
discussed in Section~\ref{Rat}, we know that 
\[
d_2(b,\cR_{n,n})=
d_2\Bigl((\bar b-\overline{ b(0)})\,,\,\cR_{n-1,n}\Bigr).
\]
Now, the Hankel operators $\Gamma_{\bar{b}}$ and 
$\Gamma_{\bar{b}-\overline{b(0)}}$ coincide and we saw in the proof of
Lemma~\ref{bab} that $\Gamma_{\bar b}$ is an isometry from 
$(\text{Ker}\,\Gamma_{\bar{b}})^\perp=\cP_{d-1}/\widetilde{q}_d$ onto 
$\text{Im}\Gamma_{\bar b}=\cP_{d-1}/\widetilde{q}_d$. Hence the
$e_j$ given by \eqref{defMW} for $1\leq j\leq d$ form an orthonormal family  of
$d\geq n+1$  singular vectors associated with the singular value 1. 
By the AAK theorem it follows that $d_{\infty}(\bar b,H^\infty_j)=1$ for $0\leq j\leq d-1$,
and since $\|e_j\|_\infty=(1-|\zeta_j|^2)^{-1/2}$, estimate \eqref{borninfB} 
follows at once from Theorem~\ref{LBH2AAK} upon choosing $v_j=e_{j+1}$ for 
$0\leq j\leq n$.

Next, if we let $w_j=(\sum_{k=0}^n e^{2i\pi kj/(n+1)}v_k)/(n+1)^{1/2}$ for
$0\leq j\leq n$, we
get another orthonormal family of $n+1$ singular vectors associated with the 
singular value 1, and clearly 
$\|w_j\|_\infty\leq (n+1)^{-1/2}\sum_{k=0}^{n+1}(1-|\zeta_k|^2)^{-1/2}$ 
for all $j$. Estimate \eqref{borninfBa} 
now follows from Theorem~\ref{LBH2AAK} again upon replacing
the previous $v_j$ by $w_j$.

If now $b$ is infinite and $k\geq0$ is the multiplicity of the zero at the 
origin,
we can write \eqref{defBlaschke} for some constant $c$ of unit modulus.
Let us define $b_m=cz^m$ if $m\leq k$ and 
\begin{equation}
\label{defBlaschkesec}
b_m(z)=cz^k\prod_{l=k+1}^m\frac{{-\bar \zeta}_l}{|\zeta_l|}\,\frac{z-\zeta_l}{1-{\bar \zeta}_lz},\qquad m>k.
\end{equation}
The sequence of Blaschke products
$\{b_m\}$ converges to $b$ pointwise on $\DD$, and since it 
is bounded it must also converge weakly to $b$ in $H^2$. Since $b_m$ and $b$ 
have norm 1, the limit of the norms is the norm of the weak limit, hence the 
convergence is actually strong in $H^2$ \cite[Theorem 3.32]{brezis}. 
Consequently  
\[\lim_{m\to\infty}\,\,d_2(b_m,\cR_{n,n})=d_2(b,\cR_{n,n}),
\]
and since estimates  \eqref{borninfB}, \eqref{borninfBa} depend only of the 
first $n+1$ zeros of $b$ they remain valid in the limit.

\hfill\boite

In view of  Corollary~\ref{corquot}, the conclusion of Theorem~\ref{LBB} 
can be sharpened upon replacing in the proof
$\|v_j\|_\infty$ and $\|w_j\|_\infty$ by $|||\Gamma_{\bar b v_j}|||$
and $|||\Gamma_{\bar b w_j}|||$. Computations become  more
involved but in any case cannot increase the left hand side
of  \eqref{borninfB} and  \eqref{borninfBa} by more than a 
factor 2. Incidentally, for 
$q_n\in\cP_N$ having all roots in $\DD$,
it seems to be an open question which $L^2$-orthonormal 
bases of 
$\cP_{n-1}/q_n$ have minimax $L^\infty$-norm.
Using such bases instead of $e_j$ in the proof of Theorem~\ref{LBB} 
may improve on  the result.

Since \eqref{Bcond} is necessary and sufficient for $\{\zeta_l\}$ to be the 
zero set of a Blaschke product, an immediate corollary to Theorem~\ref{LBB} is:
\begin{corollary}
Whenever $\alpha_n$ is a nonincreasing sequence in $(0,1]$ such that
$\Sigma_n \alpha_n<\infty$, there is a Blaschke product $b$ such that
\begin{equation}
\label{borninfBe}
\frac{\alpha_{n+1}^{1/2}}{\sqrt{n+1}}\leq 
\inf_{r\in \cR_{n,n}\cap H^{2}}\|b-r\|_2,\qquad n\in\NN,
\end{equation}
and also 
\begin{equation}
\label{borninfBae}
\left(\sum_{j=0}^{n}\frac{1}{\alpha_j^{1/2}}\right)^{-1}
\leq 
\inf_{r\in \cR_{n,n}\cap H^{2}}\|b-r\|_2,\qquad n\in\NN.
\end{equation}
\end{corollary}
\subsection{Application to Cauchy integrals on hyperbolic geodesics}
Recall that geodesic lines for the hyperbolic metric in $\DD$ 
are radii and circular arcs orthogonal to~$\TT$~\cite[ch. I]{gar}.
By definition, a hyperbolic  geodesic
segment is a compact and connected subset thereof. Alternatively, a 
hyperbolic geodesic
segment is the image of a real segment $[a,b]\subset(0,1)$ 
under an automorphism of 
the disk ({\it i.e.} a M\"obius transformation of the type 
$z\mapsto \alpha (z-z_0)/(1-\bar z_0z)$ with $|\alpha|=1$ and $z_0 \in \DD$,
in other words a Blaschke product of degree 1). Below is a nonstandard 
characterization of hyperbolic geodesic segments which is
analytic in nature. We will not use the ``if'' part but is is interesting 
in itself.

\begin{lem}
\label{caracG}
A $C^1$-smooth, closed Jordan arc $\gamma\subset\DD$ is a hyperbolic 
geodesic segment if, and only if there is a constant $C=C(\gamma)>0$ such that,
to each $g\in H^2$, there is $h\in H^2$ with $h_{|\gamma}=\bar g_{|\gamma}$
and $\|h\|_2\leq C\|g\|_2$. If $g$ is continuous on $\overline{\DD}$, so is
$h$.
\end{lem}
{\it Proof:} if $\gamma$ is hyperbolic geodesic segment,
then it is the image of a real segment under
an automorphism $\varphi$ of $\DD$ and 
$h(z)=\overline{(g\circ\varphi)(\bar z)}\circ\varphi^{-1}$
does the job. Conversely, if $\gamma$ is a $C^1$-smooth closed Jordan arc 
in $\DD$  with endpoints 
$z_1$, $z_2$ and if there exists a constant $C=C(\gamma)$ as in the statement
of the lemma, then the proof of \cite[thm. 10.1]{BS02} applies (upon
trading the geodesic arc $\bf G$  for $\gamma$ in that proof) to show that 
$\gamma$ consists exactly of non-isolated points of the cluster set,
as $n$ ranges over $\NN$, of 
poles of best approximants  to 
$((z-z_1)(z-z_2))^{-1/2}\in C(\TT)$ from $H^\infty_n$.
Because this characterization depends
only on $z_1$, $z_2$, it follows that $\gamma$ must be the geodesic arc 
joining them.

\hfill\boite

In this section, we will consider functions of the form
\begin{equation}
\label{fG}
f(z)=\frac{1}{2i\pi}\int_{G} \frac{h(\xi)}{z-\xi}d\xi
\end{equation}
where: 
\begin{itemize}
\item[(H1)] $G\subset\DD$ is a geodesic segment,
\item[(H2)] $h$ is a complex-valued function
on $G$, summable with respect to arclength, having  continuous argument
except possibly for finitely many jumps of amplitude $\pi$.
\end{itemize}
 The prototype of such a function is one which is analytic over $\DD$ 
except for two branchpoints of order strictly greater than -1. 
Indeed, by Cauchy formula,
such a function can be written as the Cauchy integral, on any 
smooth cut connecting the branchpoints, of the jump of the function across 
the cut. This jump is locally analytic and has continuous argument on the cut
(up to the branchpoints by Puiseux expansion), except at the zeros that 
the jump 
may have on this cut where  the argument has left and right limits 
which differ by $k\pi$ if $k$ is the order of the zero.
Choosing the hyperbolic geodesic cut, we get representation 
\eqref{fG}.  It may  seem artificial to 
favor the hyperbolic geodesic segment linking the branchpoints
among all possible cuts. However, this one turns out
to attract almost all poles of best rational approximants
(see \cite{BaYaSta} for this and generalizations to 
finitely many branchpoints) and also of best meromorphic approximants
(see \cite[thm. 10.1]{BS02} and Corollary~\ref{borninfw} below),
which makes it in some sense
the natural singular set of the function.

We  need additional facts from AAK theory that shed light on 
singular vectors of Hankel operators with continuous symbol.
They apply in  particular to $\Gamma_f$ when $f$ is of the form~\eqref{fG}.
\begin{itemize}
\item For $f\in C(\TT)$ and $n\geq0$, a best approximant  $g_n$ 
to $f$ from $H^\infty_n$ in $L^\infty$ uniquely exists
\cite[thm. 1.3]{AAK} \cite[ch. 4, thm 1.3]{Peller} which is
given  by
\begin{equation}
\label{formAAKapp}
g_n=\frac{{\bf P}_+(fv_n)}{v_n},\qquad f-g_n=\frac{\Gamma_f(v_n)}{v_n}
=\frac{{\bf P}_-(fv_n)}{v_n},
\end{equation}
where $v_n$ is \emph{any} singular vector of $\Gamma_f$ associated with
$s_n(\Gamma_f)$; moreover, the error function $f-g_n$ has constant 
modulus $s_n(\Gamma_f)$  a.e. on $\TT$ \cite[thm. 1.3]{AAK}
\cite[ch. 4, sec. 1, eqn. (1.12)]{Peller}. 
In particular, \eqref{formAAKapp} entails that the ratios
${\bf P}_\pm(fv_n)/v_n$ are independent of which singular vector $v_n$ 
associated with  $s_n(\Gamma_f)$ is used; this is remarkable for
if  $s_n(\Gamma_f)$ has multiplicity $\mu$, then the union of $\{0\}$ and of 
all associated singular vectors is a vector space 
of complex dimension $\mu$. 

\item When $f\in C(\TT)$, the inner factor of
a singular vector of $\Gamma_f$ is a finite Blaschke product. 
More precisely, keeping notations as in the previous item and letting
in addition $m=m(n)$ be the smallest non-negative integer such that 
$s_m(\Gamma_f)=s_n(\Gamma_f)$, the singular vector $v_n$ may be inner-outer 
factorized as
\begin{equation}
\label{iofvn}
v_n=bb_mw_n
\end{equation}
where $w_n\in H^2$ is outer and $b_m\in B_m$ is a Blaschke product of exact 
degree $m$ with zeros the poles of $g_n$ ($=g_m$), 
while $b$ is a finite Blaschke 
product whose zeros are also zeros of ${\bf P}_+(fv_n)$.
Moreover, with  $b$, $b_m$ and $w_n$ as in \eqref{iofvn},
 it holds that
\begin{equation}
\label{LE}
\Gamma_f(v_n)(z)=s_n(\Gamma_f)\,z^{-1}\overline{b_m(1/\bar z)}
\overline{j(1/\bar z)}\overline{w_n(1/\bar z)},\qquad |z|\geq 1
\end{equation}
where $j$ is a finite Blaschke product such that $jb\in B_{\mu-1}$ and
$\mu$ is the multiplicity of~$\sigma_n(\Gamma_f)$
\cite[thm. 1.2]{AAK}.
\item Assumptions and notations being as in the previous items, let
$v_n$ be a singular vector of $\Gamma_f$ associated with $s_n(\Gamma_f)$ and
\eqref{iofvn} be its inner-outer factorization. {\it We claim} that
$b_mw_n$ is also
a singular vector of $\Gamma_f$ associated with $s_n(\Gamma_f)$.
Indeed, we know from the previous item that
$g_n=b_m^{-1}h$ for some $h\in H^\infty$. 
Since $\|f-h/b_m\|_\infty=\|fb_m-h\|_\infty$,
the fact that $g_n$ is a best approximant to $f$ from $H^\infty_n$
entails that $h$ is the best approximant to 
$fb_m\in C(\TT)$ from $H^\infty$, hence
$|||\Gamma_{fb_m}|||=\|f-h/b_m\|_\infty=s_n(\Gamma_f)$ by the AAK theorem. 
Taking into account that $\Gamma_{fb_m}(u)=\Gamma_f(b_mu)$ for $u\in H^2$,
and also that $\Gamma_f^*(\Phi)={\bf P}_+(\bar f \Phi)$ for $\Phi\in H^{2,0}$,
while using that $\bar b_m \bar{H}^{2,0}\subset \bar{H}^{2,0}$ and
${\bf P}_++{\bf P}_-=Id$, we now compute
\begin{align}
\nonumber\Gamma_{fb_m}^*\Gamma_{fb_m}(bw_n)=&
\Gamma_{fb_m}^*\Gamma_f(v_n)
={\bf P}_+\Bigl(\overline{fb_m}\,\,\Gamma_{f}(v_n)  \Bigr)
={\bf P}_+\Bigl(\bar b_m{\bf P}_+\bigl(\bar f\Gamma_f(v_n)\bigr)
\Bigr)\\
\nonumber
=&{\bf P}_+\left(\bar b_m\Gamma_f^*\Gamma_f(v_n)\right)=
s^2_n(\Gamma_f){\bf P}_+\left(\bar b_mv_n\right)=
s^2_n(\Gamma_f)bw_n.
\end{align}
This shows that $bw_n$ is a maximizing vector of $\Gamma_{fb_m}$.
Next, we observe that
\begin{equation}
\label{mvdiv}
\|\Gamma_{fb_m}(bw_n)\|_2=\|{\bf P}_-\bigl(b\Gamma_{fb_m}(w_n)\bigr)\|_2
\leq \|\Gamma_{fb_m}(w_n)\|_2
\end{equation}
because multiplication by $b$ is an isometry and anti-analytic projection 
is a contraction in $L^2$. Since $\|bw_n\|_2=\|w_n\|_2$, we conclude from
\eqref{mvdiv} that $w_n$ is in turn a maximizing vector of $\Gamma_{fb_m}$ 
and that equality must hold throughout in this equation. 
In other words $b\Gamma_{fb_m}(w_n)\in H^{2,0}$, which implies easily
that  $b\Gamma_{fb_m}(w_n)=\Gamma_{fb_m}(bw_n)$. Consequently
\begin{align}
\label{simp}
\Gamma_f^*\Gamma_f(b_mw_n)&={\bf P}_+\left(
\overline{fb} \,b \Gamma_{f}(b_mw_n)\right)=
{\bf P}_+\left(\bar b{\bf P}_+(\bar f\Gamma_{f}(bb_mw_n))\right)\\
\nonumber
&=
{\bf P}_+\left(\bar b\Gamma_f^*\Gamma_f(v_m)\right)=
s_n^2(\Gamma_f){\bf P}_+\left(\bar bv_n\right)=
s_n^2(\Gamma_f)b_mw_n.
\end{align}
{\it This proves the claim}.

\end{itemize}

We now assume that $f$ has the form \eqref{fG}. Using \eqref{Cauchyproj} to 
express definition \eqref{Hankelinf} of the Hankel operator, 
then inserting \eqref{fG} and using successively Fubini's theorem and the
residue  formula, we obtain:
\begin{equation}
\label{fintvm}
\Gamma_f(v_n)(z)=\left(\frac{1}{2i\pi}\right)^2
\int_{G} h(\xi)d\xi
\int_\TT \frac{v_n(\zeta)}{(\zeta-\xi)(z-\zeta)}d\zeta=
\frac{1}{2i\pi}\int_{G}\frac{v_n(\xi)h(\xi)}{z-\xi}d\xi,\quad |z|>1.
\end{equation}
In particular $\Gamma_f(v_n)$ extends analytically from
$\overline{\CC}\setminus\overline{\DD}$ to $\overline{\CC}\setminus G$,
and Equation~\eqref{LE}  becomes
\begin{equation}
\label{repvc}
s_n(\Gamma_f)\,z^{-1}\overline{b_m(1/\bar z)}
\overline{j(1/\bar z)}\overline{w_n(1/\bar z)}=
\frac{1}{2i\pi}\int_{G}\frac{v_n(\xi)h(\xi)}{z-\xi}d\xi,\qquad |z|\geq1.
\end{equation}
Multiplying the restriction of \eqref{repvc} to $z\in\TT$ 
by $b_mj$ and  then taking anti-analytic projection again gives us
after a similar computation:
\begin{equation}
\label{repwc}
s_n(\Gamma_f)\,z^{-1}
\overline{w_n(1/\bar z)}=
\frac{1}{2i\pi}\int_{G}\frac{j(\xi)b^2_m(\xi) b(\xi) w_n(\xi)h(\xi)}{z-\xi}d\xi,\qquad |z|>1,
\end{equation}
where we took into account \eqref{iofvn}. Equation~\eqref{repwc}
entails that in turn $\check{w}_n$ ({\it cf.} \eqref{defcheck}) 
extends analytically from
$\overline{\CC}\setminus\overline{\DD}$ to $\overline{\CC}\setminus G$,
or equivalently that $w_n$ extends analytically from $\DD$ to
$\overline{\CC}\setminus\overline{G}^{-1}$, where
$\overline{G}^{-1}$ is the reflection of $G$ across $\TT$.

We can now establish a technical result which is the key for applying
Theorem~\ref{LBH2AAK}  to functions of the form \eqref{fG}.
Recall that a family of analytic functions in an open set $\Omega\subset\CC$ 
is said to be normal if it is uniformly bounded on every compact subset of 
$\Omega$. Equivalently, a normal family of analytic functions is one which is
relatively compact for the topology of locally uniform convergence in $\Omega$.

\begin{prop}
\label{normalf}
Let $f$ assume the form \eqref{fG} where hypotheses H1-H2
do hold,  and $\{v_n\}_{n\in\NN}$ be a sequence 
of singular vectors of $\Gamma_f$ such that $\|v_n\|_2=1$ for all $n$. 
Denote by $w_n$ the outer factor of~$v_n$.
Then, $\{w_n\}_{n\in\NN}$ is a normal family
in $\overline{\CC}\setminus\overline{G}^{-1}$.
\end{prop}
{\it Proof:} we already pointed out that $w_n$ is
analytic in  $\overline{\CC}\setminus\overline{G}^{-1}$.
According to Equation~\eqref{iofvn}, the inner-outer factorization of $v_n$
is of the form $v_n = bb_mw_n$, and we know 
from a previous claim ({\it cf.}~Equation~\eqref{simp}) that $b_mw_n$ is another 
singular vector of $\Gamma_f$ associated with $s_n(\Gamma_f)$ 
having the same outer factor $w_n$. Hence we can replace $v_n$ by $b_mw_n$ (in other words, we may --~and we shall~-- assume that $b\equiv 1$ and write $v_n=b_mw_n$).
For correctness, one should of course write $m(n)$ throughout, but we drop
the  dependence  of $m$ on $n$ for simplicity.

To prove that $w_n$ is bounded independently of $n$
on each compact subset of 
 $\overline{\CC}\setminus\overline{G}^{-1}$,
we parallel the argument of \cite[thm. 10.1]{BS02}.

Let $t\mapsto \alpha(t)$ parametrize $G$ with an automorphism 
$\alpha$ of $\DD$ as $t$ ranges over a real segment $[a,b]$.  
Then $t\mapsto\alpha'(t)$ has continuous 
argument. Let $\beta_1$ be a finite Blaschke product with real coefficients
vanishing precisely at the 
jumps of amplitude $\pi$ that $t\mapsto\arg h(\alpha(t))$ may have on $[a,b]$
(if $h(\alpha(t))$ is continuous we simply put $\beta\equiv1$). Then
$t\mapsto\arg (\beta_1(t)h(\alpha(t)))$ is continuous by our assumptions on
$h$. Thus, by 
Mergelyan's theorem, there is a polynomial $T$ which is real valued on 
$[a,b]$ and such that 
$|T(t)+\arg \alpha'(t)+\arg (\beta_1(t)h(\alpha(t)))|<\pi/3$ for 
$t\in[a,b]$. In invariant 
form, this means that the function 
$H=P\circ\alpha^{-1}\in H^\infty\cap C(\TT)$ is 
real valued on $G$ and moreover that 
\begin{equation}
\label{redresse}
\left|\beta(\xi)h(\xi)d\xi\right|=
\left|e^{iH(\xi)}\beta(\xi)h(\xi)d\xi\right|\leq 2\text{Re}
\left(e^{iH(\xi)}\beta(\xi)h(\xi)d\xi\right),\qquad \xi\in G,
\end{equation}
where $\beta=\beta_1\circ\alpha^{-1}$ is in turn a finite Blaschke product
which is real-valued on $G$.
Notice that $H$ and $\beta$ depend only on $f$ and not on $n$.

In another connection, 
since $w_n$ has no zero in $\DD$, it has a well-defined 
square root $w_n^{1/2}\in H^\infty$. Note, since $\|v_n\|_2=\|w_n\|_2=1$
by assumption, that $\|w_n^{1/2}\|_2=\|w_n\|_1^{1/2}\leq1$ by the Schwarz 
inequality.  Appealing to Lemma~\ref{caracG}, let
$H_n\in H^2$ take  conjugate values to $b_mw_n^{1/2}$ on $G$, with
$\|H_n\|_2\leq C\|b_mw^{1/2}_n\|_2= C\|w^{1/2}_n\|_2\leq C$. 
Note that $H_n$ is continuous on $\overline{\DD}$
since $b_m$ and $w_n^{1/2}$ are.
For $j$ as in \eqref{LE} ($j$ depends on $v_n$ but we drop this dependence), 
consider the contour integral
\begin{equation}
\label{defI}
\frac{s_n(\Gamma_f)}{2i\pi}\int_\TT e^{iH(\xi)} \,H_n(\xi)\,\beta(\xi)\,
\overline{b_m(\xi)j(\xi)}\,\frac{\overline{w_n(\xi)}}{w_n^{1/2}(\xi)}
\,\frac{d\xi}{\xi}
\end{equation}
where it should be observed that the integrand is continuous even though
$w_n$ may have zeros on $\TT$ (of course 
at such points $\bar w_n/w_n^{1/2}$ is 
understood to be $0$).
In view of 
\eqref{repvc}, this integrand 
extends analytically on $\DD\setminus G$, hence we may rewrite \eqref{defI} as
an integral over the circle $\TT_r=\{z\,:\,|z|=r\}$ where $r\in(0,1)$
is close enough to $1$ that $\TT_r$ encompasses $G$. 
Then, substituting \eqref{repvc} and using again Fubini's theorem 
and the Cauchy formula (which is permitted since $w_n^{1/2}(z)$ does not 
vanish for $|z|\leq r$), the integral \eqref{defI} transforms into
\begin{align}
\nonumber
\frac{1}{(2i\pi)^2}\int_{\TT_r} 
\left(\int_G\frac{v_n(\zeta)h(\zeta)}{\xi-\zeta}d\zeta\right)
\frac{e^{iH(\xi)} \,H_n(\xi)\beta(\xi)}{w_n^{1/2}(\xi)}
\,d\xi
=&
\frac{1}{2i\pi}\int_G \frac{v_n(\zeta)}{w_n^{1/2}(\zeta)}h(\zeta)
e^{iH(\zeta)} \,H_n(\zeta)\beta(\zeta)\,
d\zeta\\
\nonumber
&\!\!\!\!\!\!\!\!\!\!\!\!\!\!\!\!\!\!\!\!\!\!=\frac{1}{2i\pi}\int_G b_m(\zeta)w_n^{1/2}(\zeta)h(\zeta)
e^{iH(\zeta)} \,H_n(\zeta)\,\beta(\zeta)\,
d\zeta,
\end{align}
where we took \eqref{iofvn} into account. Altogether,
by the construction of $H_n$, we deduce that
\begin{equation}
\label{egdesiva}
\frac{s_n(\Gamma_f)}{2i\pi}\int_\TT e^{iH(\xi)} \,H_n(\xi)\,\beta(\xi)
\overline{b_m(\xi)j(\xi)}\,\frac{\overline{w_n(\xi)}}{w_n^{1/2}(\xi)}
\,\frac{d\xi}{\xi}=
\frac{1}{2i\pi}\int_G \bigl|b_m^2(\zeta)\,w_n(\zeta)\bigr|h(\zeta)\,
\beta(\zeta)
e^{iH(\zeta)}
d\zeta.
\end{equation}
By \eqref{redresse}, we get on the one hand that
\begin{equation}
\label{tRe}
\frac{1}{4\pi}\int_G \Bigl|b_m^2(\zeta)\,w_n(\zeta)
\,\beta(\zeta)\,h(\zeta)\Bigr|d|\zeta|
\leq
\left|\frac{1}{2i\pi}\int_G \bigl|b_m^2(\zeta)\,w_n(\zeta)\bigr|\beta(\zeta)
h(\zeta)
e^{iH(\zeta)}\right|
d\zeta.
\end{equation}
On the other hand, since $\beta,b_m,j$ are Blaschke products while 
$\|w_n^{1/2}\|_2\leq1$ and $\|H_n\|_2\leq C$, we see from
the Schwarz inequality that
\begin{equation}
\label{estinti}
\left|\frac{s_n(\Gamma_f)}{2i\pi}\int_\TT e^{iH(\xi)} \,H_n(\xi)\,\beta(\xi)\,
\overline{b_m(\xi)j(\xi)}\,\frac{\overline{w_n(\xi)}}{w_n^{1/2}(\xi)}
\,\frac{d\xi}{\xi}\right|\leq 
C \,s_n(\Gamma_f) \,\| e^{iH}\|_\infty.
\end{equation}
Therefore, in view of \eqref{tRe},  \eqref{egdesiva}, and \eqref{estinti},
we get that 
\begin{equation}
\label{bel}
\frac{1}{2\pi}\int_G \Bigl|b_m^2(\zeta)\,w_n(\zeta)
\beta(\zeta)h(\zeta)\Bigr|d|\zeta|\leq 2C \,s_n(\Gamma_f) \,\| e^{iH}\|_\infty.
\end{equation}
Now, if we multiply \eqref{repwc} (where $b\equiv1$) by $\beta$ and apply 
${\bf P}_-$ to this product,
the computation based on Fubini's theorem and Cauchy formula that led us to
\eqref{fintvm} and \eqref{repwc} yields
\begin{equation}
\label{wpb}
s_n(\Gamma_f)\,{\bf P}_-(\beta \check{w}_n)(z)=
\frac{1}{2i\pi}\int_{G}\frac{j(\xi)b^2_m(\xi) b(\xi) w_n(\xi)\beta(\xi)
h(\xi)}{z-\xi}d\xi,\qquad |z|\geq1.
\end{equation}
Equation~\eqref{wpb} entails that ${\bf P}_-(\beta \check{w}_n)$ extends 
analytically to $\overline{\CC}\setminus G$ and, 
as $|j|\leq 1$ in $\DD$ since it is a Blaschke product,
it follows from \eqref{bel} and 
\eqref{wpb}  that
\[
|{\bf P}_-(\beta \check{w}_n)(z)|\leq 2C \,\| e^{iH}\|_\infty\,
\bigl(\inf_{\zeta\in G}|z-\zeta|
\bigr)^{-1},\qquad z\in \overline{\CC}\setminus G.
\]
This proves that $|{\bf P}_-(\beta \check{w}_n)|$ is uniformly bounded
with respect to $n$ on every compact subset of $\overline{\CC}\setminus G$.
In another connection, observe from \eqref{Cauchyproj}
that  ${\bf P}_+(\beta \check{w}_n)$ is uniformly bounded with respect to $n$
on compact subsets of $\DD\setminus G$, because 
$\|\beta \check{w}_n\|_2=\|\check{w}_n\|_2=1$. Adding up, we get that
$\beta\check{w}_n$ is uniformly bounded with respect to $n$
on compact subsets of $\DD\setminus G$. Since $|\beta|$, which is a finite 
Blaschke product with all its zeros on $G$,  is bounded from
below on compact subsets of $\overline{\CC}\setminus G$,
we thus conclude that $\{\check{w}_n\}$ is normal in 
$\overline{\CC}\setminus G$. By reflection across $\TT$,
normality of $\{w_n\}$ in $\overline{\CC}\setminus\overline{G}^{-1}$ follows, 
as desired.
\hfill\boite

The following corollary to Proposition~\ref{normalf} is
worth pointing out as it shows in a rather strong sense 
that most of the poles of best $L^\infty$  meromorphic approximants to
$f$ as in \eqref{fG} asymptotically cluster to $G$.
\begin{corollary}
\label{borninfw}
\label{asfbpp}
Let $f$ assume the form \eqref{fG} where hypotheses H1-H2
do hold. Denote by $g_n$
the best approximant 
to $f$ from $H^\infty_n$ in $L^\infty$.
To each neighborhood $\mathcal{V}(G)$ of $G$, there are  
$n_0,N_0\in \NN$ such that,  if $n\geq n_0$, then
$g_n$ has at most $N_0$ poles outside $\mathcal{V}(G)$,
counting multiplicity.
\end{corollary}
{\it Proof:} we make notations as in the proof of Proposition~\ref{normalf}.
We noticed already before the latter  that $\check{w}_n$ is
analytic in  $\overline{\CC}\setminus G$.
In addition, it is clear that $ \overline{b_m(1/\bar z)}=1/b_m(z)$  
(resp. $\overline{j(1/\bar z)}=1/j(z)$) since $b_m$ (resp. $j$) is 
unimodular on $\TT$. Hence
$\overline{b_m(1/\bar z)}$  (resp. $\overline{j(1/\bar z)}$)
is meromorphic 
in $\CC$ with poles at the zeros of $b_m$ (resp. of $j$) and no zero in 
$\overline{\DD}$. Since the right hand side of   
\eqref{repvc} is analytic in $\overline{\CC}\setminus G$, we conclude that
every zero of $b_m$ (and of $j$) which does not lie on $G$
is a zero of $\check{w}_n$ with 
same or greater multiplicity. Now,  by Proposition~\ref{normalf},
every subsequence  $\check{w}_{n_k}$,
has a subsequence  $\check{w}_{n_{k_\ell}}$ converging locally 
uniformly in $\overline{\CC}\setminus G$ to some analytic function 
$\check{w}$ which is not the zero function 
because $\|\check{w}\|_2=\lim_{\ell\to\infty}\|\check{w}_{n_{k_\ell}}\|_2=1$. 
In particular
$\check{w}$ has only finitely zeros $z_1\cdots,z_N$ of respective 
multiplicities
$\mu_1,\cdots,\mu_N$ in $\DD\setminus\mathcal{V}_G$. Thus, by the Rouch\'e 
theorem, $\check{w}_{n_{k_\ell}}$ has exactly $\mu_j$ zeros in the 
neighborhood of $z_j$ for 
$\ell$ large enough, counting multiplicities, and no other zero in
in $\DD\setminus\mathcal{V}_G$. Consequently
every subsequence of $\{\check{w}_n\}$ has boundedly many zeros in  
$\DD\setminus\mathcal{V}_G$, which implies the desired conclusion as 
poles of $g_m$ which do not lie on $G$ are zeros of $\check{w}_n$  by the 
first part of the proof.

\hfill\boite

The main result of this section is the following.

\begin{theorem}
\label{asfbp}
Let $f$ assume the form \eqref{fG} where hypotheses H1-H2
do hold. Then, 
\begin{equation}
\label{comparerr}
C_1\frac{d_\infty(f,\cR_{n-1,n})}{\sqrt{n+1}}\leq C_2\frac{d_{\infty}(f,H^\infty_n)}{\sqrt{n+1}}\leq d_2(f,\cR_{n-1,n})
\end{equation}  
where $C_1,C_2$ are strictly positive constants depending on $f$ but
not on $n$.
\end{theorem}
{\it Proof:} if  $v_n$ is a singular vector of $\Gamma_f$ associated with 
$s_n(\Gamma_f)$, normalized so that $\|v_n\|_2=1$, and if $w_n$ is the outer 
factor of $v_n$, we deduce 
 from Proposition~\ref{normalf} that $\|w_n\|_\infty=\|v_n\|_\infty$ is 
bounded 
independently of $n$. Hence the second inequality in \eqref{comparerr} 
follows from Theorem~\ref{LBH2AAK}.

To prove  the first inequality, we must show that
\begin{equation}
\label{polcmer}
d_\infty(f,\cR_{n-1,n})\leq C \,d_{\infty}(f,H^\infty_n)
\end{equation}
for some constant $C$ independent of $n$. Let $g_n$
be  a best approximant 
to $f$ from $H^\infty_n$ in $L^\infty$, and write $g_n=r_n+h_n$
where $r_n\in\cR_{n-1,n}\cap \bar{H}^{\infty,0}$ while $h_n\in H^\infty$.
Note that $f-r_n\in \bar{H}^{\infty,0}$, hence $h_n={\bf P}_+(f-g_n)$.
Obviously it holds that $d_\infty(f,\cR_{n-1,n})\leq \|f-r_n\|_\infty$, 
therefore,  
it is enough to check that $\|f-r_n\|_\infty\leq C d_{\infty}(f,H^\infty_n)$
in order to establish \eqref{polcmer}. Now, by the triangle inequality,
we get that
\[\|f-r_n\|_\infty\leq \|f-g_n\|+\|h_n\|_\infty=d_{\infty}(f,H^\infty_n)+
\|{\bf P}_+(f-g_n)\|_\infty,
\]
and  we are left to prove that  
$\|{\bf P}_+(f-g_n)\|_\infty\leq C d_{\infty}(f,H^\infty_n)$.
Let $v_n$ be a singular vector of $\Gamma_f$, associated with
$s_n(\Gamma_f)$, having inner-outer factorization
$v_n=b_mw_n$, where $b_m\in B_m$ vanishes exactly at the poles of $g_n$ and
$w_n$ is outer; this is possible by a previous claim ({\it cf.} \eqref{simp}).
Here and below, we should write
for correctness $m=m(n)$, but we drop the dependence of $m$ on $n$
for simplicity.
From \eqref{formAAKapp} and \eqref{LE}, we gather that
\[{\bf P}_+(f-g_n)=s_n(\Gamma_f)\,{\bf P}_+\left(\bar b^2_m\bar j\check{w}_n/w_n
\right),\]
where we also dropped the dependence of $j$ on $v_n$,
and since $s_n(\Gamma_f)=d_{\infty}(f,H^\infty_n)$ 
it remains to establish that 
$\|{\bf P}_+\left(\bar b^2_m\bar j\check{w}_n/w_n
\right)\|_\infty$ is bounded independently of $n$. For this, 
it is enough to show that from any subsequence $n_k$ one can extract a 
subsequence  $n_{k_\ell}$ for which the property holds.
Appealing to Proposition~\ref{normalf} as
in the proof of Corollary~\ref{borninfw},   we can extract
from $\{w_{n_k}\}$ a subsequence $\{w_{n_{k_\ell}}\}$
converging locally uniformly to some 
$w$, analytic in 
$\overline{\CC}\setminus\overline{G}^{-1}$, which is not the
zero function.
Let  $w$ have $N$ zeros lying on $\TT$, say
$z_1,\cdots,z_N$, where multiplicities are accounted by repetition
and it is understood if $N=0$ that $\{z_j\}$ is the 
empty set. Pick $\varepsilon>0$ small enough that
the circle $\TT_{1+\varepsilon}$ does not meet
$\overline{G}^{-1}$ and 
$w$ has no other zeros than $z_1,\cdots,z_N$ in the corona 
$\mathcal{C}_\varepsilon=\{z:\,1\leq |z|\leq1+\varepsilon\}$.
When $\ell$ is large enough, by the Rouch\'e theorem,
$w_{n_{k_\ell}}$ has exactly $N$ zeros $z_{1,\ell},\cdots, z_{N,\ell}$ in 
$\mathcal{C}_\varepsilon$, counting multiplicities with repetition,
and $\{z_{j,\ell}\}$ converges to $\{z_j\}$ as a set when
$\ell\to+\infty$
(recall that $w_{n_{k_\ell}}$ is outer hence has no zero in $\DD$).
We label the $z_{j,\ell}$ so that, say $z_{j,\ell}\in\TT$ for $1\leq j\leq s_\ell$ and
$z_{j,\ell}\notin\TT$ for $s_\ell+1\leq j\leq N$.
Define
\[P_{N,\ell}(\xi)=\Pi_{j=1}^N(\xi-z_{j,\ell}),\qquad
Q_{N-s_\ell,\ell}(\xi)=\Pi_{j=s_\ell+1}^N(\xi-z_{j,\ell}),
\]
and let us write $w_{n_{k_\ell}}(\xi)=u_\ell(\xi)\,P_{N,\ell}(\xi)$
where  $u_\ell(\xi)$ is 
analytic in $\overline{\CC}\setminus\overline{G}^{-1}$ and zero-free in
$\mathcal{C}_\varepsilon$. By the maximum principle,
$u_{j,\ell}$  converges to $w(z)/\Pi_{j=1}^N(z-z_j)$ locally uniformly
in  $\overline{\CC}\setminus\overline{G}^{-1}$.
Clearly,
\begin{equation}
\label{expquow}
\frac{\check{w}_n(\xi)}{w_n(\xi)}=
\Pi_{l=1}^{s_\ell}\left(-\bar z_{l,\ell}\right)\,
\xi^{-N}\,
\frac{\check{u}_{j,\ell}(\xi)}{u_{j,\ell}(\xi)}\,
\frac{\widetilde{Q}_{N-s_\ell,\ell}(\xi)}{Q_{N-s_\ell,\ell}(\xi)},
\end{equation}
where we observe that 
$b_{N-s_\ell}=\Pi_{l=1}^{s_\ell}\left(-\bar z_{l,\ell}\right)\,\widetilde{Q}_{N-s_\ell,\ell}/Q_{N-s_\ell,\ell}$ lies in
$B_{N-s_\ell}$ and that
$\check{u}_{j,\ell}/u_{j,\ell}$ is continuous and
bounded independently of $\ell$ on
$\mathcal{C}_\varepsilon$ as well as analytic in the interior of $\mathcal{C}_\varepsilon$. Now, put $\beta(\xi)=\xi$ and let us write
\begin{equation}
\label{decompP+}
{\bf P}_+\left(\bar b^2_m\bar j\check{w}_m/w_m
\right)=b_{N-s_\ell}{\bf P}_+\left(\bar\beta^N\bar b^2_m\bar j\check{u_\ell}/u_\ell
\right)+
{\bf P}_+\left(b_{N-s_\ell}{\bf P}_-\left(\bar\beta^N\bar b^2_m\bar j\check{u}_\ell/u_\ell
\right)\right).
\end{equation}
Recalling that  $\overline{b_m(1/\bar z)}=1/b_m(z)$ and 
$\overline{j(1/\bar z)}=1/j(z)$,  we deduce 
from 
\eqref{Cauchyproj}
\begin{equation}
\label{fintP+err}
{\bf P}_+\left(\bar\beta^N\bar b^2_m\bar j\check{u_\ell}/u_\ell\right)(z)
=\frac{1}{2i\pi}\int_\TT
\frac{1}{b^2_m(\xi) j(\xi)}
\frac{\check{u}_\ell(\xi)}
{u_\ell(\xi)(\xi-z)}\,\frac{d\xi}{\xi^N},
\qquad |z|<1,
\end{equation}
and by Cauchy's theorem we can deform the contour of integration to
$\TT_{1+\varepsilon}$ without changing the value of the integral:
\begin{equation}
\label{fintP+errb}
{\bf P}_+\left(\bar\beta^N\bar b^2_m\bar j\check{u_\ell}/u_\ell\right)(z)
=\frac{1}{2i\pi}\int_{\TT_{1+\varepsilon}}
\frac{1}{b^2_m(\xi) j(\xi)}
\frac{\check{u}_\ell(\xi)}
{u_\ell(\xi)(\xi-z)}\,\frac{d\xi}{\xi^N},
\qquad |z|<1.
\end{equation}
The integral in the right hand side of \eqref{fintP+errb}
is now bounded in modulus,  independently of $\ell$ and $z\in\DD$, 
because $|b_m|\geq1$ and $| j|\geq1$ 
on  $\overline{\CC}\setminus\DD$, while $\varepsilon\leq |\xi-z|$ and 
$\check{u}_\ell/u_\ell$ is uniformly bounded on $\TT_{1+\varepsilon}$. 
Thus, the first summand in the right hand side of \eqref{decompP+}
is bounded in $L^\infty$, independently of $\ell$, because
$b_{N-s_\ell} \in B_{N-s_\ell}$. To see that the second summand is also 
bounded, we put 
$\Psi={\bf P}_-\left(\bar\beta^N\bar b^2_m\bar j\check{u}_\ell/u_\ell
\right)$ and  we notice that $\|\Psi\|_\infty$ is bounded independently 
of $\ell$ because $\bar\beta^N\bar b^2_m\bar j\check{u}_\ell/u_\ell$ is
unimodular on $\TT$ and we just saw from \eqref{fintP+errb} that
$\|{\bf P}_+\left(\bar\beta^N\bar b^2_m\bar j\check{u}_\ell/u_\ell
\right)\|_\infty$ is bounded independently of $\ell$.
Next, we observe that this second summand is
${\bf P}_+(b_{N-s_\ell}\Psi)$ and that it lies in
$\cR_{N-s_\ell-1,N-s_\ell}$ because
for $\xi\in\TT$ we have:
\[\frac{Q_{N-s_\ell}(\xi)}{\xi^{N-s_\ell}}{\bf P}_+(b_{N-s_\ell}\Psi)(\xi)=
\overline{Q_{N-s_\ell}(1/\bar\xi)}\Psi(\xi)-
\frac{Q_{N-s_\ell}(\xi)}{\xi^{N-s_\ell}}{\bf P}_-(b_{N-s_\ell}\Psi)(\xi),
\]
and both summands on the right  lie in $\bar{H}^{2,0}$ whence
$Q_{N-s_\ell}{\bf P}_+(b_{N-s_\ell}\Psi)\in\cP_{N-s_\ell-1}$.
We now appeal to Grigoryan's theorem, saying that if 
${\bf P}_-(\Phi)\in\mathcal{R}_{d-1,d}$ then 
$\|{\bf P}_-(\Phi)\|_\infty\leq cd\|\Phi\|_\infty$ for some absolute constant 
$c$, see  \cite[eqn. (6.1)]{Peller}. As 
$\check{{\bf P}}_+(\Phi)={\bf P}_-(\check{\Phi})$ for 
any function $\Phi$, it implies since the check operation preserves
$\cR_{N-s_\ell-1,N-s_\ell}$ and the $L^\infty$ norm that
\[\|{\bf P}_+(b_{N-s_\ell}\Psi)\|_\infty
=\|{\bf P}_-(\bar b_{N-s_\ell}\check{\Psi})\|_\infty
\leq c(N-s_\ell)\|\bar b_{N-s_\ell}\check{\Psi}\|_\infty=c(N-s_\ell)
\|\Psi\|_\infty.
\]
This achieves the proof.

\hfill\boite

The authors conjecture that Theorem~\ref{asfbp} carries over to Cauchy
integrals of the form \eqref{fG} where $G$ is a so-called symmetric contour 
for the Green potential in $\DD$ ({\it cf.} \cite[thm 1]{St85b} for details),
the prototype of which is an analytic 
function with finitely many branchpoints of order greater than $-1$
in the disk.  For such functions, more generally
even if branchpoints have arbitrary order, it was proved in \cite{GRakh87} that
$\lim_{n\to\infty}d_\infty(f,\cR_{n-1,n})^{1/n}=\exp\{-2/C\}$ where $C$ is the
condenser capacity of the pair $(\TT,G)$. The same $n$-th root estimate holds 
for $d_2(f,\cR_{n-1,n})$ \cite[cor. 8]{BaYaSta}, and more generally for
the distance from $f$ to $\cR_{n-1,n}$ in $L^p$  when 
$1\leq p\leq\infty$ \cite{Totik}.
Inequality \eqref{comparerr} compares  $d_2(f,\cR_{n-1,n})$ and
$d_\infty(f,\cR_{n-1,n})$ in a much stronger sense, but still one may
wonder if the factor $1/\sqrt{n+1}$ is really needed. 
In the special case of Markov functions,  {\it i.e.} of Cauchy integrals of 
positive densities on a segment, the results in \cite{BSW} show that
this factor  is in fact superfluous.

\section{Linearized errors}
\label{LES}
Given $f\in \bar{H}^{2,0}$,
 $p_{n-1}\in \cP_{n-1}$, $q_n\in \cP_n$ and a (complex)
weight function $w\in H^\infty$, the
\emph{linearized error} associated with $p$, $q$ and $w$ in problem
RAB($n$) is
\begin{equation}
\label{deflin}
\mathcal{L}(f,p_{n-1},q_n,w):=(q_nf-p_{n-1})w.
\end{equation}
It is formally obtained from the error $f-p_{n-1}/q_n$ by chasing 
denominator $q_n$ and multiplying by the weight.
In applied sciences, problem RAB($n$) and weighted variants thereof 
are of great 
importance to model time series  as well as to identify  linear 
dynamical systems\footnote{
For continuous time systems, rational approximation is
performed on the imaginary axis rather than the circle.
This, is equivalent to the present setting thanks to the 
isometry $f\mapsto
\sqrt{2}f((z+1)/(z-1))/(z-1)$ mapping $H^2$ onto the Hardy space of 
$\{{\rm Re}z>0\}$  while preserving rationality and the degree.}, 
{\it e.g.} in modal analysis of mechanical structures or
in frequency analysis of microwave devices
\cite{HD,HerVan,Ljung,Regalia,OSM}. The importance of the $L^2$ norm in this 
context stems from its statistical interpretation as a variance.
Because  RAB($n$) (or equivalently MA($n$)) is a difficult non convex problem,
several  approaches to system identification in engineering  have been
based on linearization.
Most popular in this connection are two closely related heuristics,
namely the Steiglitz-McBride method \cite{SMB,Regalia,PAMM} and
the vector fitting method \cite{GuSe,DHD}.
These are iterative procedures, first   choosing $w=1/\pi_n$ where
$\pi_n\in \cP_n$ is monic with no root on $\TT$,
then minimizing $\|\mathcal{L}(f,p_{n-1},q_n,w)\|_2$
with respect to 
$p_{n-1}$ as well as $q_n$, the latter being normalized so as to be monic\footnote{What we describe here is the Steiglitz-McBride method, although
we should mention that the criterion used is often a discretized version of 
$\|\mathcal{L}(f,p_{n-1},q_n,1/\pi_n)\|_2$ obtained from pointwise values
on $\TT$. The vector fitting method is essentially a rewriting of the 
Steiglitz-McBride procedure where rational functions are parametrized in
pole-residue form.
}  
(which yields a convex problem). Subsequently, one 
replaces $w$ by $1/q^{o}_n$, where $q_n^{o}$ is the optimal $q_n$, 
and repeats the previous steps until some fixed point is reached. 
Such procedures are prompted by the easy observation that if 
$f\in R_{n-1,n}$,
then the value of the problem is zero from the first iteration already.
Accordingly, convergence was studied in a classical
stochastic setting for system identification, where 
$f\in R_{n-1,n}\cap \bar{H}^{2,0}$ is perturbed by white noise 
(the noise is then constitutive of the model), but
such heuristics do not converge in 
general when $f\notin R_{n-1,n}$ 
\cite{RMA,SLAA,PAMM}.

Our purpose here is not to discuss these techniques, nor to compare them with
dedicated optimization algorithms \cite{FO,OSM},  but rather to stress a  
link between the value of RAB($n$) and the minimization 
of linearized errors.

We consider weights of the form
$w=1/\pi_n$ where $\pi_n\in \cP_n$ is a polynomial having no root on $\TT$.
When minimizing the linearized error, Theorem~\ref{singlust}
suggests a specific normalization for $q_n$: let us define
\begin{equation}
\label{normalization}
P_{\pi_n}:=\{q_n\in \cP_n:\ \|q_n/\pi_n\|_\infty=1\}.
\end{equation}
Then, the following result holds.
\begin{theorem}
\label{thmpi}
Let $f\in \bar{H}^{2,0}$ and $\pi_n\in \cP_n$, with 
$\mathcal{Z}(\pi_n)\cap\TT=\emptyset$.
Then 
\begin{equation}
\label{borninfpi}
 d_2(f,\cR_{n-1,n})=d_2(f,H^2_n)
 \geq
 \min_{\stackrel{q_n\in P_{\pi_n}}{p_{n-1}\in \cP_{n-1} }}
\|\mathcal{L}(f,p_{n-1},q_{n,}1/\pi_n)\|_2.
\end{equation}
\end{theorem}
{\it Proof:}
we may assume without loss of generality that $\mathcal{Z}(\pi_n)\subset
\overline{\CC}\setminus\DD$, for otherwise we can replace
every linear factor $(z-a)$ of $\pi_n$ for which  $a\in\DD$ by the linear
factor $(1-\bar az)$ which has reflected zero across $\TT$. This leaves
$|\pi_n|$ unchanged on $\TT$, and consequently does not affect the 
minimization of $\|\mathcal{L}(f,p_{n-1},q_n,1/\pi_n)\|_2$.

Now, arguing as we did to obtain \eqref{minp}, 
we find that $\pi_n{\bf P}_+(fq_n/\pi_n)$ is a 
polynomial of degree at most $n-1$. 
If we write $p_{n-1}(q_n,\pi_n)$ for this polynomial and take into account that
$\mathcal{Z}(\pi_n)\cap\overline{\DD}=\emptyset$, we see from 
Parseval's theorem 
that for fixed $q_n\in\cP_n$ the criterion
$\|\mathcal{L}(f,p_{n-1},q_{n},1/\pi_n)\|$ gets minimized 
precisely when $p_{n-1}=p_{n-1}(q_n,\pi_n)$, so that
\begin{equation}
\label{explinerH}
\min_{p_{n-1}\in\cP_{n-1}} \|\mathcal{L}(f,p_{n-1},q_n,1/\pi_n)\|_2=\|{\bf P}_-
\left(f\frac{q_n}{\pi_n}\right)\|_2=\|A_f(q_n/\pi_n)\|_2.
\end{equation}
Let
\[K_{\pi_n}:=\{q_n/\pi_n:\ q_n\in P_{\pi_n}\}.\]
Since $\pi_n$ has no zeros on $\overline{\DD}$, it holds that
$K_{\pi_n}\subset  \mathcal{S}^\infty$, the unit sphere of $H^\infty$.
Identifying $\cP_n$ with $\CC^{n+1}\sim\RR^{2n+2}$ by taking 
coefficients  as coordinates, we see that  $K_{\pi_n}$
is homeomorphic to the Euclidean sphere $\sph^{2n+1}$ via 
the map $q_n/\pi_n\mapsto q_n/\|q_n\|_2$ which is 
odd. Therefore $K_{\pi_n}$ is a compact subset of $\mathcal{S}^\infty$
of genus $2n+2$, and by  \eqref{critsinglust}:
\begin{equation}
\label{minsxi}
 d_2(f, R_{n-1,n})=
d_2(f,H^2_n)\geq \min_{q_n/\pi_n \in K_{\pi_n}} \|A_f(q_n/\pi_n)\|_2
\end{equation}
which is \eqref{borninfpi} in view of \eqref{explinerH}.

\hfill\boite

It follows easily from a compactness argument that the minimum  in
the right hand side of \eqref{borninfpi} is attained. 
However, it not {\it a priori}
obvious  how to compute it for
$P_{\pi_n}$ is not convex. Numerically, this issue can be approached 
as follows.  First, we assume without loss of generality that
$\pi_n$ has no roots in $\overline{\DD}$, so that
\eqref{explinerH} holds ({\it cf.} proof of Theorem~\ref{thmpi}).
Next, for $\xi\in\TT$, let
 \[P_{\pi_n,\xi}:=\{q_n\in \cP_n:\ \|q_n/\pi_n\|_\infty=1,\ q_n(\xi)=\pi_n(\xi)\}.\]
Observe that $P_{\pi_n,\xi}$ is never empty  when $\pi_n$ has no zero 
on $\TT$. Indeed, for small $\varepsilon>0$, it holds that
$|\pi_n(e^{i\theta})|^2-|\varepsilon(e^{i\theta}-\xi)|^2\geq0$
hence, by Fej\`er-Riesz factorization
(see Lemma~\ref{FR} to come), there is a polynomial $q_n$ with
$|q_n|\leq|\pi_n|$ on $\TT$ and $|q_n(\xi)|=|\pi_n(\xi)|$. Thus,
$q_n\pi_n(\xi)/q_n(\xi)$ lies in  $P_{\pi_n,\xi}$.
Clearly 
\begin{equation}
\label{decompK}
K_{\pi_n}=\cup_{\xi,\,\zeta\in\TT} \,\,\zeta \,P_{\pi_n,\xi},
\end{equation}
and multiplying $q_n$ by $\zeta\in\TT$ cannot change the value of
$\|A_f(q_n/\pi_n)\|_2$. Therefore it holds that
\begin{equation}
\label{minpsi}
 \min_{\stackrel{q_n\in P_{\pi_n}}{p_{n-1}\in \cP_{n-1} }}
\|\mathcal{L}(f,p_{n-1},q_n,1/\pi_n)\|_2=
  \min_{\xi\in\TT} \psi(\xi)
\end{equation}
where the function $\psi(\xi)$ is given by ({\it cf.} \eqref{explinerH})
\begin{equation}
\label{defpsi}
\psi(\xi)=\min_{ \stackrel{q_n\in P_{\pi_n,\xi}}{p_{n-1}\in \cP_{n-1} }}
 \|\mathcal{L}(f,p_{n-1},q_n,1/\pi_n)\|_2=
\min_{ q_n\in P_{\pi_n,\xi}}\|A_f(q_n/\pi_n)\|_2.
\end{equation}
Note that $\psi(\xi)$ can be computed as the solution of a convex problem 
for each $\xi$, because $ P_{\pi_n,\xi}$ is a convex set and 
$\|A_f(q_n/\pi_n)\|_2$ a quadratic criterion. Granted this ability
to evaluate $\psi$ pointwise, we discuss below how to
numerically estimate the minimum in \eqref{minpsi}.

Clearly $\psi$ is the zero function when
$f\in\cR_{n-1,n}$, for if $f=p/q$ with $\text{deg}\,q\leq n$ we may pick 
$q_n=q$ as minimizer in \eqref{defpsi}. The next lemma describes 
this minimizer in greater detail when $f\notin\cR_{n-2,n-1}$.
\begin{lem}
\label{cgm}
Let $f\in \bar{H}^{2,0}$ and $\pi_n\in \cP_n$, with 
$\mathcal{Z}(\pi_n)\cap\overline{\DD}=\emptyset$.
If $f\notin\cR_{n-2,n-1}$, then
the minimizing $q_n$ in \eqref{defpsi} is unique, has 
all its roots in $\overline{\DD}$, and exact degree $n$.
\end{lem}
{\it Proof:}  assume that
$q_{n,1}$ and  $q_{n,2}$  are distinct minimizers, that is,
$q_{n,1}, q_{n,2}\in P_{\pi_n,\xi}$ and
$\|A_f(q_{n,1}/\pi_n)\|_2=\|A_f(q_{n,2}/\pi_n)\|_2=\psi(\xi)$. Put 
$q_{n,3}=(q_{n,1}+q_{n,2})/2\in P_{\pi_n,\xi}$. 
By strict convexity of the $L^2$ 
norm, we get $A_f(q_{n,1}/\pi_n)=A_f(q_{n,2}/\pi_n)$
otherwise we would have that $\|A_f(q_{n,3}/\pi_n)\|_2<\psi(\xi)$
which is absurd. Set $q=q_{n,1}-q_{n,2}\in \cP_n$.
Then $A_f(q/\pi_n)=0$ implying by
definition of $A_f$ that $fq/\pi_n\in H^2$.
A fortiori then  $fq\in H^2$, and since $f\in \bar{H}^{2,0}$ we must
have that $fq$ is a polynomial of degree at most
$n-1$, say $p$. Thus, $f=p/q$, and as $q$ has a root on $\TT$
(namely $\xi$) the latter must be cancelled by a corresponding root of $p$.
Altogether $f\in\cR_{n-2,n-1}$, thereby
showing the uniqueness part of the lemma.
Let now $q_{n,\xi}$ be the unique minimizer and $b$ 
a Blaschke product with poles in 
$\mathcal{Z}(q_{n,\xi})\cap\overline{\CC}\setminus\overline{\DD}$.
Then $bq_{n,\xi}\in\cP_{n}$ has same
modulus as $q_n$ on $\TT$, hence there is $\zeta\in\TT$ such that
$\zeta bq_{n,\xi}\in P_{\pi_n,\xi}$. Since $\zeta b$ is a Blaschke product,
reasoning as in \eqref{mvdiv} yields 
$\|A_f(\zeta b q_{n,\xi}/\pi_n)\|_2\leq \|A_f(q_{n,\xi}/\pi_n)\|_2$
so that
$\zeta b q_{n,\xi}$ is in turn a minimizer, hence is equal to
$q_{n,\xi}$ by the uniqueness part just proved. Now,
$q_{n,\xi}\not\equiv0$ since $q_{n,\xi}(\xi)=\pi_n(\xi)\neq0$, therefore
$\zeta b=1$. Thus, $b$ must be a constant, that is to say
there cannot be a zero of $q_{n,\xi}$ outside $\overline{\DD}$. 
Finally, assume that $\text{deg}\,q_{n,\xi}<n$. Then $q_{n,\xi}(z)z\bar\xi$ 
lies in  $P_{\pi_n,\xi}$ and,
since $z\bar\xi$ is a Blaschke product, it follows as before that 
$q_{n,\xi}z\bar\xi$ is a minimizer, hence it must be equal to $q_{n,\xi}$ by
uniqueness. This  contradiction achieves the proof.

\hfill\boite

We need a continuity property of the  Fej\`er-Riesz factorization
that we could not ferret out in 
the literature.
Write $\mathcal{T}_n$ for the space of trigonometric polynomials of degree at 
most $n$, {\it i.e.} sums of the form $\sum_{|k|\leq n} a_ke^{ik\theta}$.
For fixed $n$, $\cP_n$ and $\mathcal{T}_n$ have a natural topology
induced by any norm.
\begin{lem}
\label{FR}
To each nonzero $T\in\mathcal{T}_{n}$ such  that $T\geq0$ on $\TT$, one can associate
continuously a unique polynomial $q\in\cP_n$ having no zero in 
${\DD}$ and such that $|q(e^{i\theta})|^2=T(e^{i\theta})$ with $q(0)>0$.
\end{lem}
{\it Proof:} 
 let
$\mathcal{T}_n^+\subset \mathcal{T}_n$ be the closed subset of 
trigonometric  polynomials which are non-negative on $\TT$.
For $T\in\mathcal{T}_n^+$, 
existence of $q\in \cP_n$ such that $|q|^2=T$ on $\TT$
is a classical result known after Fej\`er and Riesz
\cite[sec. 53]{RN}. 
Since $|z-a|=|1-z\bar a|$ for $z\in\TT$,
clearly $q$ may be chosen zero free in $\DD$ if  $T\not\equiv0$.
Then, $q_{|\DD}$ is outer  in $H^\infty$, for it has 
no zero and it extends analytically across
$\TT$ \cite[ch. II, thms. 6.2 \& 6.3]{gar}. Thus,
formula \eqref{defext}
shows that $q$ is uniquely defined by $\log|q|=\log T/2$, therefore also by 
$T$. Moreover,
each coefficient of $q$ is a 
continuous function of $\log T\in L^1$, because the $k$-th coefficient is just
the derivative 
$q^{(k)}(0)/k!$ and we may differentiate  \eqref{defext} 
under the integral sign. To achieve  the proof, we
establish that $T\mapsto\log T$ 
is continuous from $\mathcal{T}^+\setminus\{0\}$ into $L^1$.

First, {\it we claim} that $T\mapsto\|\log T\|_1$ is continuous
from $\mathcal{T}_n^+\setminus\{0\}$ into $\RR$.  
To see this, it is enough to show
that if  $T^{\{k\}}$ tends to
$T$ in $\mathcal{T}_n^+\setminus\{0\}$ as $k\to\infty$, then 
$\|\log T^{\{k_\ell\}}\|_1$ tends to $\|\log T\|_1$ for some subsequence
$T^{\{k_\ell\}}$. By the first part of the proof, we can write 
$T^{\{k\}}=|q^{\{k\}}|^2$ with
\[q^{\{k\}}(z)=q^{\{k\}}(0)\,\Pi_{l=1}^n(1-z \,a^{\{k\}}_l),\qquad a^{\{k\}}\in\overline{\DD},
\]
where multiplicities are counted by repetition and the ordering of the roots
for each $k$ is arbitrary. Note that $|q^{\{k\}}(0)|$ is 
bounded, since by the Schwarz inequality: 
\[|q^{\{k\}}(0)|=\left|\frac{1}{2\pi}
\int_0^{2\pi} q^{\{k\}}(e^{i\theta})d\theta\right|\leq
\|q^{\{k\}}\|_2=\left\|T^{\{k\}}\right\|_1^{1/2}.\]
Therefore, there is a subsequence $q^{\{k_\ell\}}$ such that
$q^{\{k_\ell\}}(0)$ converges to $c\in\CC$ and 
$a^{\{k_\ell\}}_j$ converges to $a_j\in\overline{\DD}$  for each 
$j\in\{1,\cdots,n\}$. If we let
\[q(z)=c\,\Pi_{j=1}^n\,(1-z \,a_j),\]
then clearly $|q^{\{k_\ell\}}|^2$ converges to $|q|^2$ almost everywhere on 
$\TT$, so that necessarily $|q|^2=T$. 
In particular, we have that $c\neq0$ otherwise $T$ would be identically zero,
a contradiction.
Now, since $\log$ turns products into sums, 
we are left to show that if $b^{\{k\}}\to b$ in $\overline{\DD}$, then
\begin{equation}
\label{limlogfact}
\lim_{k\to+\infty}\int_0^{2\pi}\left|\log|1-e^{i\theta}b^{\{k\}}|\right|\,
d\theta=
\int_0^{2\pi}\left|\log|1-e^{i\theta}b|\right|\,d\theta.
\end{equation}
When $|b|<1$ relation \eqref{limlogfact} is obvious. If $|b|=1$,
we may assume by rotational symmetry that $b=1$ and  $b^{\{k\}}\in[0,1]$,
in which case \eqref{limlogfact} follows by dominated convergence from the 
observation 
that $|1-b^{\{k\}}e^{i\theta}|\geq |\sin\theta|$ for $|\theta|\leq\pi/2$.
{\it This proves  the claim.}

Now, if  $T^{\{k\}}$ tends to
$T$ in $\mathcal{T}_n^+\setminus\{0\}$, it is plain that
$\log T^{\{k\}}$ converges to $\log T$ almost everywhere on $\TT$
and by the previous claim 
the $L^1$-norm of the limit is the limit of the $L^1$-norms.
Thus, the desired $L^1$-convergence of $\log T^{\{k\}}$ to $\log T$  
follows from Egoroff's theorem 
\cite[ch.3, ex.17]{Rudin}.

\hfill\boite

With the help of Lemmas~\ref{cgm} and \ref{FR}, we  now prove that 
$\psi$ is continuous:
\begin{lem}
\label{psicont}
Let $f\in \bar{H}^{2,0}$ and $\pi_n\in \cP_n$, with 
$\mathcal{Z}(\pi_n)\cap\overline{\DD}=\emptyset$.
Then the map $\psi$ defined by \eqref{defpsi} is continuous on $\TT$.
\end{lem}
{\it Proof:} if $f\in \mathcal{R}_{n-1,n}$, we mentioned
that $\psi\equiv0$ already.
Otherwise, dwelling on Lemma~\ref{cgm}, 
let $q_{n,\xi}$ indicate the unique minimizer in the last term of
\eqref{defpsi}. By definition of $P_{\pi_n,\xi}$ we have that
$|q_{n,\xi}|\leq |\pi_n|$ on $\TT$ and that $q_{n,\xi}(\xi)=\pi_n(\xi)$, 
in particular $q_{n,\xi}$ is
bounded independently of $\xi$. Thus, from any convergent sequence
$\xi_k\to\xi$ on $\TT$, 
we can extract a subsequence $\xi_{k_\ell}$ for which 
$q_{n,\xi_{k_\ell}}$ converges uniformly to
some $q\in\cP_{n}$, and passing to the limit we see that
$q\in P_{\pi_n,\xi}$. Given $\varepsilon>0$, we can
pick the sequence $\xi_k$ so that
\[\lim_{k\to+\infty}\psi(\xi_k)=l\leq\liminf_{\zeta\to\xi} \psi(\zeta)+\varepsilon,\]
and by continuity of $q_n\mapsto\|A_f(q_n/\pi_n)\|_2$ from $\cP_n$ into $\RR$
we get that
\begin{equation}
\label{scipsi}
\liminf_{\zeta\to\xi}\psi(\zeta)+\varepsilon
\geq l=\lim_{\ell\to\infty}\|A_f(q_{n,\xi_{k_\ell}}/\pi_n)\|_2=
\|A_f(q/\pi_n)\|_2\geq
\psi(\xi).
\end{equation}
Since $\varepsilon>0$ was arbitrary, we conclude  that
$\psi$ is lower semi-continuous. To see that $\psi$
is in fact continuous, it is enough to establish the {\it following claim}:
to each $\xi\in\TT$ and  
$\varepsilon>0$, there is $\eta>0$ such that 
$|\xi-\zeta|<\eta$ implies existence of
$q_\zeta\in P_{\pi_n,\zeta}$ with 
$\|q_{n,\xi}-q_\zeta\|_\infty<\varepsilon$. 
Indeed, if the claim holds, 
we get from 
\eqref{defpsi} that when $|\xi-\zeta|<\eta$: 
\[\psi(\zeta)\leq\|A_f(q_\zeta/\pi_n)\|_2\leq 
\|A_f(q_{n,\xi}/\pi_n)\|_2+\|A_f((q_\zeta-q_{n,\xi})/\pi_n)\|_2
\leq\psi(\xi)+\varepsilon\|f/\pi_n\|_2,\]
and since $\varepsilon$ was arbitrary we conclude that
$\limsup_{\zeta\to\xi}\psi(\zeta)\leq \psi(\xi)$ whence $\psi$
is indeed continuous in view of \eqref{scipsi}.

To establish the claim,  observe from Lemma~\ref{FR} since
$|\pi_n|^2-|q_{n,\xi}|^2$ is a non-negative trigonometric polynomial 
of degree at most $n$ on $\TT$ that 
\begin{equation}
\label{appFR}
|q_{n,\xi}(e^{i\theta})|^2+|\kappa_{n,\xi}(e^{i\theta})|^2=|\pi_n(e^{i\theta})|^2
\end{equation}
where $\kappa_{n,\xi}\in\cP_n$ has no root in $\DD$ and is uniquely defined
by \eqref{appFR} together with  the normalization
$\kappa_{n,\xi}(0)>0$.
As $\pi_n(\xi)=q_{n,\xi}(\xi)$, we can write 
$\kappa_{n,\xi}(z)=(z-\xi)Q_{n-1}(z)$, and for $\zeta\in\TT$ we set more 
generally: 
\begin{equation}
\label{defkappanzeta}
\kappa_{n,\zeta}(z)=\frac{(z-\zeta)Q_{n-1}(z)}{\lambda_\zeta},
\qquad \lambda_\zeta=\sup_{z\in\TT}|(z-\zeta)Q_{n-1}(z)/\pi_n(z)|.
\end{equation}
Clearly  $|\kappa_{n,\zeta}|^2\leq |\pi_n|^2$ on $\TT$ 
so that, by Lemma~\ref{FR}, there is a unique $P_{n,\zeta}\in \cP_n$ having no 
root  in $\DD$ and meeting $P_{n,\zeta}(0)>0$ such that 
\begin{equation}
\label{defPnzeta}
|P_{n,\zeta}(e^{i\theta})|^2+|\kappa_{n,\zeta}(e^{i\theta})|^2=
|\pi_n(e^{i\theta})|^2.
\end{equation}
Since $q_{n,\xi}$ has all its roots in $\DD$ and exact degree $n$ by 
Lemma~\ref{cgm}, it follows from \eqref{appFR} and the uniqueness part of
Lemma~\ref{FR} that
$P_{n,\xi}=\bar c\widetilde{q}_{n,\xi}$ where $c\in\TT$ is such that 
$cq_{n,\xi}$ 
has positive leading coefficient. Moreover, it is easily checked from
\eqref{defkappanzeta} that $|\kappa_{n\zeta}|^2$ is arbitrary close to
$|\kappa_{n,\xi}|^2$ in $\mathcal{T}_n$ if $|\xi-\zeta|$ is sufficiently 
small. Therefore, from \eqref{defPnzeta}
and the continuity property asserted by Lemma~\ref{FR},
we deduce that $P_{n,\zeta}$ is arbitrary close to 
$\bar c\widetilde{q}_{n,\xi}$ in $\cP_n$ 
if $|\xi-\zeta|$ is sufficiently  small. Consequently
$Q_\zeta=\bar c\widetilde{P}_{n,\zeta}$ is
arbitrary close to $q_{n,\xi}$  when $|\xi-\zeta|$ is small enough.
Now, by \eqref{defPnzeta} and the definition of $Q_\zeta$,
it holds that $|Q_\zeta|=|P_{n,\zeta}|\leq|\pi_n|$ on $\TT$,
and also that $|Q_\zeta(\zeta)|=|P_{n,\zeta}(\zeta)|=|\pi_n(\zeta)|$ because
$\kappa_{n,\zeta}(\zeta)=0$ by construction.
Hence, $q_\zeta=(\pi_n(\zeta)/Q_\zeta(\zeta))Q_\zeta$ lies in 
$P_{\pi_n,\zeta}$, and since $Q_\zeta(\zeta)\to q_{n,\xi}(\xi)=\pi_n(\xi)$
as $\zeta\to\xi$ we have that $(\pi_n(\zeta)/Q_\zeta(\zeta))\to1$
when $\zeta\to\xi$. Thus, just like $Q_\zeta$, the polynomial
$q_\zeta$  is arbitrary close 
to $q_{n,\xi}$  when $|\xi-\zeta|$ is small enough, which {\it proves the claim}.

\hfill\boite

To estimate the right hand side \eqref{minpsi}, it remains to minimize
$\psi(\xi)$ over $\xi\in\TT$, which can be numerically performed
 by dichotomy because $\TT$ is compact and 1-dimensional while
$\psi$ is continuous by Lemma~\ref{psicont}.

A natural question is whether the lower bound \eqref{borninfpi} can be 
sharp. The answer is no except in the trivial case where
$f\in\mathcal{R}_{n-1,n}$:

\begin{prop}
Assumptions and notations as in Theorem~\ref{thmpi}, it holds if
$f\notin\mathcal{R}_{n-1,n}$ that
\begin{equation}
\label{borninfps}
 d_2(f, \mathcal{R}_{n-1,n})>
\,\,\min_{\stackrel{q_n\in P_{\pi_n}}{p_{n-1}\in \cP_{n-1} }}
\|\mathcal{L}(f,p_{n-1},q_n,1/\pi_n)\|_2
\end{equation}
\end{prop}
{\it Proof:}  as in the proof of Theorem~\ref{thmpi}, 
we may assume that $\mathcal{Z}(\pi_n)\subset
\overline{\CC}\setminus\DD$ and then, by \eqref{explinerH}, 
the right hand side of \eqref{borninfps} is equal to
$\min_{q_n\in P_{\pi_n}}\,\|A_f(q_n/\pi_n)\|_2$.
 Let $q_{n,0}$ be a minimizer of the latter, and
$b_{n,1}=q_{n,1}/\widetilde{q}_{n,1}$ a minimizing Blaschke product in
\eqref{minHankB}.
Multiplying $q_{n,1}$ and $\widetilde{q}_{n,1}$ by a real constant,
we may assume $q_{n,1}\in P_{\pi_{n}}$.
We can also multiply $q_{n,0}$ by a unimodular constant so that, 
using \eqref{decompK}, there is $\xi_0\in\TT$
for which $q_{n,0}\in P_{\pi_n,\xi_0}$.  
Now, if \eqref{borninfps} is an equality, we get
by definition  of 
$q_{n,0}$, $q_{n,1}$ that
\begin{align}
\nonumber
 d_2(f, \mathcal{R}_{n-1,n})=&
\|{\bf P}_-(fq_{n,0}/\pi_{n})\|_2\leq
\|{\bf P}_-(fq_{n,1}/\pi_{n})\|_2\\
\label{intercomp}
=& \|{\bf P}_-\Bigl(f(q_{n,1}/\widetilde{q}_{n,1})(\widetilde{q}_{n,1}/\pi_{n})
\Bigr)\|_2
\leq \|{\bf P}_-(fq_{n,1}/\widetilde{q}_{n,1})\|_2\\
\nonumber
=&
 d_2(f, \mathcal{R}_{n-1,n})
\end{align}
where we used in the second inequality that 
$\widetilde{q}_{n,1}/\pi_{n}\in\mathcal{S}^\infty$. 
Consequently equality holds throughout \eqref{intercomp},
implying in  particular that $|\widetilde{q}_{n,1}|=|\pi_{n}|$ on $\TT$.
Thus, as both 
polynomial have no root in $\DD$ their ratio is a unimodular
constant, and  renormalizing $q_{n,1}$ if necessary we may assume that 
$\widetilde{q}_{n,1}=\pi_{n}$. Then 
$c=\widetilde{q}_{n,1}(\xi_0)/q_{n,1}(\xi_0)$ is a unimodular constant such 
that $cq_{n,1}\in P_{\pi_n,\xi_0}$, and by the uniqueness part in 
Lemma~\ref{cgm} we see that the first inequality in \eqref{intercomp} can be 
an equality only if $q_{n,0}=cq_{n,1}$.
Altogether $q_{n,0}/\pi_n=cb_{n,1}$ is in turn an optimal Blaschke product
in \eqref{minHankB}. This optimality entails 
that \cite[thm. 8.2]{BS02} 
\[A_f^*A_f(q_{n,0}/\pi_n)
={\bf P}_+
\Bigl(|A_f(q_{n,0}/\pi_n)|^2q_{n,0}/\pi_n\Bigr).
\]
Consequently, letting 
$\langle u,\, v\rangle_{L^2}=\text{Re}\langle u,v\rangle$, 
it holds for all $q_n\in P_{\pi_n,\xi_0}$ that
\begin{align}
\nonumber 
\langle A_f(q_{n,0}/\pi_n)\,&,\,A_f(q_{n,0}/\pi_n)-A_f(q_n/\pi_n) 
\rangle_{L^2}=
\langle A^*_fA_f(q_{n,0}/\pi_n)\,,\,q_{n,0}/\pi_n-q_n/\pi_n \rangle_{L^2}\\
\label{chainecrit}
=&\langle {\bf P}_+
\Bigl(|A_f(q_{n,0}/\pi_n)|^2q_{n,0}/\pi_n\Bigr)\,,\,
q_{n,0}/\pi_n-q_n/\pi_n\rangle_{L^2}\\
\nonumber
=&\langle|A_f(q_{n,0}/\pi_n)|^2q_{n,0}/\pi_n\,,\,
q_{n,0}/\pi_n-q_n/\pi_n \rangle_{L^2}\\
\nonumber
=&\langle|A_f(q_{n,0}/\pi_n)|^2\,,\,
1-q_n/q_{n,0} \rangle_{L^2},
\end{align}
where we used in the third line that $(q_{n,0}-q_n)/\pi_n\in H^2$
to get rid of ${\bf P}_+$ and in the last line that
$\overline{q_{n,0}/\pi_n}=\pi_n/q_{n,0}$ on $\TT$.

On the one hand, we see that \eqref{chainecrit} is non-negative for
all $q_n\in P_{\pi_n,\xi_0}$ because $|q_n/q_{n,0}|=|q_n/\pi_n|\leq1$ on
$\TT$. In fact, it can be made strictly positive: indeed, 
$A_f(q_{n,0}/\pi_n)$ is not the zero function for $f\notin\mathcal{R}_{n-1,n}$,
and as discussed before 
\eqref{decompK} one can pick $q_n\in P_{\pi_n,\xi_0}$ such that
$|q_n/q_{n,0}|(\xi)<1$ for $\xi\neq\xi_0$ on $\TT$. 
On the other hand, the fact that $q_{n,0}$ is a minimizer in the
convex problem \eqref{defpsi} (where $\xi$ is set to $\xi_0$) implies
that \eqref{chainecrit} is nonpositive 
for all $q_n\in P_{\pi_n,\xi_0}$
\cite[prop. 5.23]{Demengel}, a contradiction which concludes the proof.

\hfill\boite

\section{Numerical results}
\label{sec:numericalresults}
In order to study how effective the bounds given by Theorem~\ref{LBH2AAK}, Corollary~\ref{corquot} and Theorem~\ref{thmpi}, we wrote a prototype 
implementation in each case and ran it on a few examples. 
We report in this section the results obtained 
on the following set of functions $f \in \bar{H}^{2,0}$. For each 
of them, we consider the problem of best $H^2$ approximation by a 
rational function in ${\cR}_{n-1, n}$ with $n=4$.
\begin{itemize}
\item \textbf{Example 1:} $f : z \mapsto \log( (10z-9)/(10z+9) )$.
\item \textbf{Example 2:} $f$ is a rational function of degree 5.
\item \textbf{Example 3:} $f$ is a rational function with 20 poles that have been randomly and uniformly chosen inside the unit disc.
\item \textbf{Example 4:} $f$ is a rational function with 20 poles that have been randomly and uniformly chosen inside the disc of radius $0.2$ centered at the origin.
\item \textbf{Example 5:} $f$ is a rational function with 20 poles that have been randomly and uniformly chosen inside the annulus of radii $18/20$ and $19/20$ centered at the origin.
\item \textbf{Example 6:} $f$ is fairly close to a rational function of degree 4. Namely, we chose a rational function $g$ of degree~$4$ and then we obtained $f$ by perturbing each Fourier coefficients of $g$ with a small noise of relative error bounded by $0.01$.
\item \textbf{Example 7:} $f : z \mapsto \exp(-\text{i}/(z-0.9\text{i}))-1$.
\end{itemize}
These examples have been chosen so as to exhibit
different kind of singularities inside the disc, which may or may not
be close to the 
unit circle, in order to cover various situations. 

Before discussing the results, let us say
a few words on the implementation. We are using Matlab~R2011b. For each example, $f$ is actually approximated by a truncated Fourier series 
$\widehat{f}$. The order of truncation is chosen so as to ensure that $f$ and $\widehat{f}$ agree to at least $40$ bits on the unit circle. The bounds of Theorem~\ref{LBH2AAK} and Corollary~\ref{corquot} are computed as described in Section~\ref{sec:applicationRational}: indeed, since $\widehat{f}$ is a truncated Fourier series, it can be written $\widehat{f} = p/q$ where $p \in {\cP}_{N-1}$ and $q = z^N$. It turns out that, when $q$ is a power of $z$, 
the construction in Section~\ref{sec:applicationRational} gets simpler 
since $\widetilde{q}=1$, whence Bezout relation is just 
$a\widetilde{q}+bq = 1$ with $a=1$ and $b=0$. 
We handled all examples using this 
technique, even though in examples $2$ to $5$ the number $N$ 
becomes quite large (up to $1500$) and it would have been 
more efficient (but would also have required more implementation) 
to forget about $\widehat{f}$  and to
apply the construction in Section~\ref{sec:applicationRational}
to the original $f$. Anyway, even for fairly large $N$, 
we obtained our results in a few seconds on an Intel Xeon at 2.67GHz, with 4GB of memory. When computing the bound of Theorem~\ref{LBH2AAK}, the norms $\|v_j\|_{\infty}$ are estimated by sampling $v_j$ at $8000$ points evenly distributed on the unit circle.

Regarding the bound of Theorem~\ref{thmpi}, its computation
reduces to finding the minimum of $\psi(\xi)$, for $\xi \in \TT$ (see \eqref{minpsi} and~\eqref{defpsi}) as explained in Section~\ref{LES}. For a given $\xi$, we evaluate $\psi(\xi)$ by solving a convex optimization problem. For this purpose we use CVX, a package for specifying and solving convex programs within Matlab~\cite{cvxweb, cvxpaper}. More precisely, we precompute $A_f(z^{j-1}/\pi_n)$ ($j=1\dots n+1$) in the Fourier basis, \emph{i.e.} we compute a $N\times (n+1)$ matrix $M = (m_{ij})$ such that $A_f(z^{j-1}/\pi_n) = \sum_{i=1}^N m_{ij}\, z^{-i}$. Therefore, if $q_n = \sum_{j=1}^{n+1} a_j\,z^{j-1}$, we have that
\begin{equation*}
  \|A_f(q_n/\pi_n)\|_2^2 = \left\|M
   \left(\begin{array}{c}
    a_1\\ \vdots \\ a_{n+1}
  \end{array}\right)
\right\|_2^2.
\end{equation*}
Since $N$ is much larger than $n$, it is convenient to compute a 
decomposition $M=QR$, where $Q$ is orthogonal and $R$ is upper-triangular. Since $\|Mv\|_2 = \|Rv\|_2$ for all vectors $v$ and only the first $n+1$ 
rows of $R$ are non-zero, we end up handling a $(n+1)\times(n+1)$ 
matrix instead of $M$. Now,
\begin{equation*}
  \psi(\xi) = \min_{q_n \in P_{\pi_n, \xi}} \left\|R
  \left(\begin{array}{c}
    a_1\\ \vdots \\ a_{n+1}
  \end{array}\right)
\right\|_2.
\end{equation*}
The set $P_{\pi_n,\, \xi}$ is convex, but described by infinitely many constraints. Therefore, we consider a set $\TT_1$ of $50$ points regularly spaced on the unit circle and we define
\begin{equation*}
  P_{\pi_n,\, \xi}^{(1)} = \{q_n = \sum_{j=1}^{n+1} a_j\,z^{j-1} : q_n(\xi) = \pi_n(\xi),\,\text{and } \forall \zeta \in \TT_1,\, |q_n(\zeta)| \le \pi_n(\zeta)\}.
\end{equation*}
We first use CVX to solve our minimization problem subject to $q_n \in P_{\pi_n,\,\xi}^{(1)}$. This gives an optimal polynomial $q_n^{(1)}$. Next, we construct a set $\TT_2$ by adding to $\TT_1$ the points of $\TT$ where $|q_n^{(1)}/\pi_n|$ reaches a local maximum. We then use CVX to solve our minimization 
problem subject to $q_n \in P_{\pi_n,\,\xi}^{(2)}$, where
\begin{equation*}
  P_{\pi_n,\, \xi}^{(2)} = \{q_n = \sum_{j=1}^{n+1} a_j\,z^{j-1} : q_n(\xi) = \pi_n(\xi),\, \forall \zeta \in \TT_2,\, |q_n(\zeta)| \le \pi_n(\zeta)\}.
\end{equation*}
This gives a new optimal polynomial $q_n^{(2)}$, and we repeat the process
until we reach a step $k$ where $\max_{\zeta \in \TT} |q_n^{(k)}(\zeta)/\pi_n(\zeta)|-1$ falls below the level of numerical errors produced by~CVX.

It is worth pointing out that, although sufficient in most cases
to get an idea of the numerical value of $\psi(\xi)$,
this procedure yields no certified estimate of the bound 
in Theorem~\ref{thmpi}. Actually, CVX is a user-friendly 
generic software, able to tackle many types of convex optimization 
problems with a powerful syntax. However, it offers little control on the 
difference between the true mathematical solution and the numerical estimate thereof.
Moreover, when too many constraints enter the game, 
it quickly  yields no solution at all. In addition,
it is probably much slower than would be  a dedicated tool 
to solve that particular convex problem. 
The point we want here to make is that accurately
estimatig the bound in Theorem~\ref{thmpi} ({\it i.e.} aiming at more than 
a prototypical illustration of the content of the paper)  
requires further work.

The numerical results proper are reported in Table~\ref{results} where the second column is $\frac{M_n(f)}{\sqrt{n+1}}$ (lower bound given by Theorem~\ref{LBH2AAK}), the third column is $\frac{Q_n(f)}{\sqrt{n+1}}$ (lower bound given by Corollary~\ref{corquot}), the fourth and fifth columns are $\min_{q_n\in P_{\pi_n},\, p_{n-1}\in \cP_{n-1}}\|\mathcal{L}(f,p_{n-1},q_{n,}1/\pi_n)\|_2$ (lower bound given by Theorem~\ref{thmpi}) for two different choices of $\pi_n$. 

For an appraisal of the sharpness of our results,
we also ran RARL2\footnote{\url{http://www-sop.inria.fr/apics/RARL2/rarl2.html}}, a software tool that tries to compute a solution to problem RAB($n$). RARL2 looks for  local minima of the criterion in a fairly systematic way, and returns the best approximant it could find. As a consequence, it gives an upper bound for  the value of problem RAB($n$) which
is likely to be tight and therefore interesting to compare with our lower bounds. The error $\|f-r\|_2$ generated by the candidate best approximant
$r$ computed  by RARL2 is reported in the last column.

The bound given by Theorem~\ref{thmpi} has the advantage of allowing the user to choose a weight $\pi_n$ which offers extra-flexibility to try to
improve the estimate. Yet, it is not obvious how to pick $\pi_n$ in general. The simplest choice is $\pi_n \equiv 1$ (reported in the fourth column of the table). Another, appealing possibility is to put $\pi_n = \widetilde{q}_n^{\;*}$ where $q_n^{\;*}$ is the denominator of the rational function computed by RARL2 (since $q_n^{\;*}$ has all its poles inside the disc, $\pi_n$ has all its poles outside, as required). The corresponding results
are reported in the fifth column of the table.

\begin{table}[htp]
  \centering
  \begin{tabular}{|c|c|c|c|c|c|}
    \hline
    \multirow{2}{*}{\small{Example}} & \multirow{2}{*}{\small{Bound of Th.~\ref{LBH2AAK}}} & \multirow{2}{*}{\small{Bound of Corollary~\ref{corquot}}} & \small{Bound of Th.~\ref{thmpi}} & \small{Bound of Th.~\ref{thmpi}} & \multirow{2}{*}{\small{RARL2}}\\
           &                 &                   & \small{with $\pi_n = 1$} & \small{with $\pi_n = \widetilde{q}_n^{\;*}$} & \\ \hline
    1      &   2.884744e-3   &    2.887532e-3    &      4.04e-3     & 10.8e-3  & 11.5e-3 \\
    2      &   7.731880e-2   &    7.732037e-2    &      12.4e-2     & 24.3e-2  & 24.72e-2\\
    3      &    2.459346     &     2.470149      &       2.286      &  0.258   & 16.6907 \\
    4      &    1.234503     &     1.234861      &       1.94       &   1.8    & 6.5721  \\
    5      &    47.26312     &     47.30424      &       2.14       &   N/A    & 178.3152 \\
    6      &   2.894007e-3   &    2.894380e-3    &      9.62e-3     & 12.46e-3 & 12.5e-3 \\
    7      &   1.780707e-4   &    1.782276e-4    &     0.7977e-4    & 6.3409e-4& 6.3742e-4 \\ \hline
  \end{tabular}
  \caption{Numerical results}
  \label{results}
\end{table}

As can be seen from the table, the refinement of Corollary~\ref{corquot} with respect to Theorem~\ref{LBH2AAK} is almost negligible on all examples. The bound given by 
Theorem~\ref{thmpi} with $\pi_n = 1$ is better than Theorem~\ref{LBH2AAK} and Corollary~\ref{corquot} in 4 cases out of 7,
but not overly so. Considering that the computation time is generally much longer for this  bound, this improvement can be considered as rather expensive. In contrast, when $\pi_n = \widetilde{q}_n^{\;*}$, the bound of Theorem~\ref{thmpi} becomes fairly sharp in cases 1,2,6, and 7, which is encouraging. This better bound comes at the cost of a longer computation time though, mostly because the convex optimization problems involved with this choice of $\pi_n$ seem more difficult to solve. In  Example~5, for instance, we were
not able to obtain reliable results from CVX.
However, it also appears
that choosing $\pi_n = \widetilde{q}_n^{\;*}$ is not \emph{always} best,
and it would be quite interesting  to further understand
which  $\pi_n$ are efficient in this respect.

\bibliographystyle{plain}
\bibliography{biblio}
\end{document}